\documentclass[a4paper, abstract=on]{scrartcl}

\usepackage[english]{babel}
\usepackage[utf8]{inputenc}
\usepackage{geometry}
\usepackage{graphicx, latexsym}
\usepackage{amsmath,amssymb,amsthm,mathdots,bm}
\allowdisplaybreaks{}
\usepackage{algorithm}
\usepackage{algorithmic}

\usepackage{framed}
\usepackage{multirow}
\usepackage{epstopdf}
\usepackage{psfrag}
\usepackage{hyperref}
\usepackage{pdfpages}
\usepackage{xspace}
\usepackage{microtype}
\usepackage{lmodern}
\usepackage{paralist}

\usepackage{setspace}
\usepackage{multicol}
\usepackage{array}
\usepackage{csquotes}
\usepackage{scrhack}
\usepackage{booktabs}
\usepackage{breakcites}
\usepackage{caption}
\usepackage{listings}
\usepackage[mathscr]{euscript}

\usepackage{cite}
\newtheoremstyle{break}     
{0.5em}         
{0.5em}         
{\itshape}    
{}              
{\bfseries}     
{.}             
{\newline}      
{}              
\theoremstyle{break}

\newtheorem{theorem}{Theorem}[section]

\newtheorem{lemma}{Lemma}[section]

\newtheorem{remark}{Remark}[section]

\renewcommand{\proofname}{Proof}

\newcommand{\norm}[1]{|| #1 ||}

\newcommand{\summe}[2]{\sum\limits_{#1}^{#2}}

\newcommand{\integraltwo}[2]{\int\limits_{#1}^{#2}}

\newcommand{\Rfield}[0]{\mathbb{R}}
\newcommand{\Cfield}[0]{\mathbb{C}}
\newcommand{\Rmat}[2]{\mathbb{R}^{#1 \times{#2}}}
\newcommand{\Cmat}[2]{\mathbb{C}^{#1 \times{#2}}}

\newcommand{\Cm}{\mathbb{C}^{-}}

\newcommand{\Cfun}[1]{\mathcal{C}{(#1)}}

\newcommand{\vareps}[0]{\varepsilon}

\newcommand{\spektrum}[1]{\Lambda({#1})}

\DeclareMathOperator{\rank}{rank}

\DeclareMathOperator{\range}{range}

\newcommand{\kron}{\otimes}

\newcommand{\dx}[1]{\mathrm{d}#1}

\newcommand{\rz}[1]{%
    \ensuremath{
        \mathord{
            \mathrm{
                #1%
            }%
        }%
    }%
}

\usepackage[software,hardware]{mymacros}

\usepackage{tikz}
\usepackage{pgfplotstable}
\usepackage{pgfplots}
\pgfplotsset{compat=1.13}
\usetikzlibrary{mindmap,shadows}
\usetikzlibrary{external}
\pgfplotsset{filter discard warning=false}

\makeatletter
\newcommand{\skipitems}[1]{%
    \addtocounter{\@enumctr}{#1}%
}
\makeatother

\newif\ifdraft
\ifdraft
    
\else
    
\fi

\author{Maximilian Behr \and Peter Benner \and Jan Heiland}
\title{Solution Formulas for Differential Sylvester and Lyapunov Equations}
\date{\today}
\begin{document}
\maketitle
\begin{abstract}
    \noindent
    The differential Sylvester equation and its symmetric version, the differential Lyapunov equation, appear in different fields of applied mathematics like control theory, system theory, and model order reduction. The few available straight-forward numerical approaches if applied to large-scale systems come with prohibitively large storage requirements. This shortage motivates us to summarize and explore existing solution formulas for these equations. We develop a unifying approach based on the spectral theorem for normal operators like the Sylvester operator $\mathcal S(X)=AX+XB$ and derive a formula for its norm using an induced operator norm based on the spectrum of $A$ and $B$.
In view of numerical approximations, we propose an algorithm that identifies a suitable Krylov subspace using Taylor series and use a projection to approximate the solution.
Numerical results for large-scale differential Lyapunov equations are presented in the last sections.

\end{abstract}

\newcommand{\BDFonecolor}{red}
\newcommand{\BDFtwocolor}{orange}
\newcommand{\BDFthreecolor}{magenta}
\newcommand{\BDFfourcolor}{brown}
\newcommand{\BDFfivecolor}{blue}
\newcommand{\BDFsixcolor}{cyan}

\newcommand{\ODEFourFiveColor}{red}
\newcommand{\ODETwoThreeColor}{orange}
\newcommand{\ODEOneOneThreeColor}{magenta}
\newcommand{\ODEOneFiveSColor}{black}
\newcommand{\ODETwoThreeSColor}{brown}
\newcommand{\ODETwoThreeTColor}{blue}
\newcommand{\ODETwoThreeTBColor}{cyan}

\tableofcontents
\section{Introduction}\label{sec_introduction}
For coefficient matrices $A\in \Cmat{n}{n}$ and $B \in  \Cmat{m}{m}$, an inhomogeneity $C\in \Cmat{n}{m}$, and an initial value $D\in \Cmat nm$, we consider the differential matrix equation
\begin{equation}
    \begin{aligned}\label{eq:introsylvestereq}
        \dot{X}(t) & = AX(t) +  X(t)B +  C, \\
        X(t_0)     & = D,
    \end{aligned}
\end{equation}
and provide formulas for the solution $X$ with $X(t)\in \Cmat{n}{m}$.
Equation~\eqref{eq:introsylvestereq} is commonly known as differential Sylvester equation (as opposed to the algebraic \emph{Sylvester equation} $AX+BX+C=0$).
In the symmetric case that $B=A^T$, equation~\eqref{eq:introsylvestereq} and its algebraic counterpart is called differential (algebraic) Lyapunov equation.
In what follows, we will occasionally abbreviate Sylvester or Lyapunov equations by SLE\@.

In particular the differential Lyapunov equation is a useful tool for stability analysis
and controller design for linear time-varying systems~\cite{AmaAAetal14}.
Equilibrium points of the differential Lyapunov equation, namely solutions of the algebraic Lyapunov equation, are used to construct quadratic Lyapunov functions for asymptotically stable linear time-invariant systems~\cite[Thm.~7.4.7]{PolW98}.
The controllability and observability problem for linear time-varying systems is strongly connected
to the solution of the differential Lyapunov equation~\cite[Ch.~13-14]{Bro70},~\cite[Ch.~3-4]{KnoK85}.
Other important applications lie in model order reduction~\cite{morAnt05} or in optimal control of linear time-invariant systems on finite time horizons~\cite{Loc01}.
Despite its importance, there have been but a few efforts to solve the differential Sylvester / Lyapunov or Riccati equation numerically, see
~\cite{LBenM04,LBenM13,Hei14,LKoeLS16,LLanMS15,LMen07,Sti15,LLan17}.
These algorithms are usually based on applying a time discretization and solving the resulting algebraic equations. Thus, even if the algebraic SLE are solved efficiently, the storage needed for the discrete solution at the time instances makes these direct approaches infeasible for large scale settings.
Recently, Krylov subspace based methods were proposed in~\cite{GueHJetal17,LHacJ17,LHacJ18,KosM17}.

In an attempt to overcome this shortage, we revisit known solution formulas, develop alternative solution representations, and discuss their suitability for numerical approximations.
We start with deriving a spectral decomposition for the Sylvester operator $\mathcal S$ which allows functional calculus. We obtain formulas for the operator norm $\|\mathcal S\|$ as well as for $\mathcal S^{-1}$ and $e^{\mathcal S}$. This recovers previously known solution formulas. It will turn out that, in terms of efficiency, this solution representation is not well suited for approximation in general but, in special cases, allows for the construction or computation of exact solutions.

As a step towards efficient solution approximation, we use Taylor series expansions to identify suitable Krylov subspaces. For the differential Lyapunov equation with stable coefficient matrices
and symmetric low-rank factored inhomogeneity, it is well-known that the solution of the algebraic Lyapunov equation spans a Krylov
subspace under these assumptions.
We split the solution of the differential Lyapunov equation in a constant and time dependent part, where
the constant part is the solution of an algebraic Lyapunov equation.
We approximate the time dependent part using the subspace spanned by the solution
of the algebraic Lyapunov equation.
The resulting algorithm overcomes the essential problems with storage consumption.
Numerical results are presented in the Section~\ref{sec_numerical_res} and Appendices~\ref{sec_appendix_numerical_res_proj} and~\ref{sec_appendix_bdf}.

\section{Preliminaries}\label{sec_prelim}
In this section, we introduce the considered equations and the Sylvester operator, set the notation, and recall basic results.

The spectrum of a matrix $A \in \Cmat{n}{n}$ is denoted by $\spektrum{A}$. 
Generally, the spectrum is a subset of $\Cfield$. A matrix is called stable if its spectrum is contained in the left open half plane $\Cm$. 
The Frobenius inner product on $\Cmat{n}{m}$ is given by $\langle A, B \rangle_F :=
    \sum\limits_{i=1}^{n} \sum\limits_{j=1}^{m} A_{i,j} \overline{B_{i,j}}$.
The Hadamard product is $ A\odot B = {\left( A_{i,j} \cdot B_{i,j} \right)}_{\substack{i=1,\ldots,n\\j=1,\ldots,m}} \in \Cmat{n}{m}$ for $A,B \in \Cmat{n}{m}$.
The Hermitian transpose, transpose, conjugate are denoted by $A^H$, $\overline{A}$, $A^T$, respectively.
Further, we will refer to the Kronecker delta: $\delta_{{ij}}={\begin{cases}1&{\text{if }}i=j, \\ 0&{\text{if }}i\neq j \end{cases}}$
and the Kronecker product: $A \kron B={\left( A_{i,j} \cdot B \right)}_{\substack{i=1,\ldots,n\\j=1,\ldots,m}}$.
The identity matrix in $\Cmat{d}{d}$ is denoted by $E_{d,d}$.

For further reference we start with recalling a general well known result on the solution which applies to the more general  case of the SLE (where the coefficient matrices may depend on time):

\begin{theorem}[Existence and Uniqueness of Solutions, {\cite[Thm.~1.1.1., Thm.~1.1.3., Thm.~1.1.5]{AboFIJ03}}]\label{thm:existence_timevariying}
    Let $I\subseteq \Rfield$ an open interval with $t_0\in I$, $A\in  \Cfun{I,\Cmat{n}{n}}, B\in \Cfun{I,\Cmat{m}{m}}, C\in \Cfun{I,\Cmat{n}{m}}$
    and $D\in \Cmat{n}{m}$.
    The differential Sylvester equation
    \begin{align*}
        \dot{X}(t) & = A(t)X(t) +  X(t)B(t) +  C(t), \\
        X(t_0)     & = D,
    \end{align*}
    has the unique solution
    $X(t) = \Phi_{A}(t,t_0)D \Phi_{B^{H}}{(t,t_0)}^H + \int\limits_{t_0}^{t} \Phi_{A}(t,s)C(s)\Phi_{B^{H}}{(t,s)}^H \dx s$. \\
    $\Phi_{A}(t,t_0)$ and $\Phi_{B^{H}}(t,t_0)$ are the unique state-transition matrices with respect to $t_0\in I$ defined by
    \begin{align*}
        \dot{\Phi}_{A}(t,t_0)   & := \frac{\partial}{\partial t} \Phi_{A}(t,t_0)=A(t)\Phi_{A}(t,t_0),         \\
        \Phi_{A}(t_0,t_0)       & \phantom{:}= E_{n,n}.                                                       \\
        \dot{\Phi}_{B^H}(t,t_0) & := \frac{\partial}{\partial t} \Phi_{B^H}(t,t_0)={B(t)}^H\Phi_{B^H}(t,t_0), \\
        \Phi_{B^H}(t_0,t_0)     & \phantom{:}= E_{m,m}.
    \end{align*}
\end{theorem}

The specification to the autonomous case with constant coefficients is straight forward by
simply replacing the state transition matrices with the matrix exponentials.

\begin{theorem}\label{thm_autonom_variation_of_const}
    Let $I\subseteq \mathbb R$ an open interval with $t_0\in I$, $A\in \Cmat{n}{n}, B\in \Cmat{m}{m}, C\in \Cfun{I,\Cmat{n}{m}}$
    and $D\in \Cmat{m}{n}$.
    The differential Sylvester equation
    \begin{align*}
        \dot{X}(t) & = AX(t) +  X(t)B +  C(t), \\
        X(t_0)     & = D,
    \end{align*}
    has the unique solution
    \begin{align}\label{eq:varofconstsAXXB}
        X(t) = e^{A(t-t_0)}D e^{B(t-t_0)} + \int\limits_{t_0}^{t} e^{A(t-s)}C(s)e^{B(t-s)}\dx s.
    \end{align}
\end{theorem}

Next, we state basic properties of the Sylvester operator, which, for given $A\in \Cmat nn$ and $B\in \Cmat mm$, is defined through its action on an $X\in\Cmat nm$:
\begin{equation}\label{eq:sylvop}
    \mathcal S(X)=A X  + XB.
\end{equation}

The Sylvester operator $\mathcal S$ has been thoroughly studied in~\cite{ByeN87,KonMP00,KonGMetal03,Ste01}.
Among others, it has been shown that the eigenvalues and eigenvectors of the Sylvester operator $\mathcal S$ can be expressed in terms of the eigenvalues and eigenvectors of $A$ and $B$,
cf.~\cite[Rem.~1.1.2.]{AboFIJ03},~\cite[Ch.~6.1.1]{morAnt05},~\cite{Koh08}.

In view of rewriting the solution formula~\eqref{eq:varofconstsAXXB}, we state the following lemma:
\begin{lemma}[Sylvester Operator $\mathcal S$]\label{lem:sylvopprops}
    For the Sylvester operator $\mathcal S\colon \Cmat{n}{m} \to \Cmat{n}{m}$ and its partial realizations $\mathcal H, \mathcal V\colon \Cmat{n}{m} \to \Cmat{n}{m}$, $\mathcal H(X)=AX$ and $\mathcal V(X) = XB$, it holds that:
    \begin{itemize}
        \item $\mathcal S = \mathcal H + \mathcal V$ and $\mathcal H \mathcal V = \mathcal V \mathcal H$,
        \item $e^{t\mathcal S}= e^{t\mathcal H}e^{t\mathcal V}$ for all $ t\in \Rfield$,
    \end{itemize}
    for any $A\in \Cmat nn $ and $B\in \Cmat mm$.
\end{lemma}
\begin{proof}
    The first claim can be confirmed by direct computations.
    The second claim is a standard result for commuting linear operators.
\end{proof}
By Lemma~\ref{lem:sylvopprops}, formula~\eqref{eq:varofconstsAXXB} rewrites as
\begin{align}
    X(t) & =  e^{tA}De^{tB} + \integraltwo{t_0}{t} e^{(t-s)A}C(s)e^{(t-s)B}\dx {s}
    =  e^{(t-t_0)\mathcal H}e^{(t-t_0)\mathcal V}(D) + \int\limits_{t_0}^{t} e^{(t-s)\mathcal H} e^{(t-s)\mathcal V}(C(s))\dx {s}            \notag \\
         & =  e^{(t-t_0)\mathcal S}(D) + \int\limits_{t_0}^{t} e^{(t-s)\mathcal S}(C(s))\dx {s}. \label{eq:variatconstSylv}
\end{align}

\section{Spectral Decomposition of the Sylvester Operator}\label{sec_spectral_sylvester}

In this section we show that the Sylvester operator $\mathcal S$, as defined in~\eqref{eq:sylvop}, is a normal operator if $A$ and $B$ are diagonalizable and if a suitably chosen inner product on a  Hilbert space is considered. The inner product depends on the decomposition of $A$ and $B$.
Nevertheless, this approach will enable us to apply the spectral theorem and to derive a  solution formula for the differential and algebraic SLE\@.
This resembles the formulas of~\cite[Ch.~4.1.1]{GajQ95},~\cite{Jam68},~\cite{LRom69}.
Those results were obtained by inserting the spectral decomposition into the SLE and
by applying suitable algebraic manipulations and using the unrolled Kronecker representation of the SLE\@.
Our strategy is to decompose the operator $\mathcal S$ first and then using functional calculus to obtain formulas for $e^{\mathcal S}$ and $\mathcal S^{-1}$.
The eigenspaces of $\mathcal S$ can be constructed from the eigenspaces of $A$ and $B$.
The choice of the inner product ensures that the eigenvectors are orthonormal
and $\mathcal S$ becomes a normal operator.

\begin{lemma}[Inner product for $\mathcal S$]\label{lemma:innerproduct}
    Let $A\in \Cmat{n}{n}, B \in \Cmat{m}{m}$  be diagonalizable and $C,D\in \Cmat{n}{m}$.
    Furthermore let $A=UD_{A}U^{-1}$ and $B^H=VD_{B^H}V^{-1}$ be the spectral decompositions of $A$ and $B^H$ with
    $U\in \Cmat{n}{n},V\in \Cmat{m}{m}$,
    $D_A$ and $D_{B^{H}}$ diagonal matrices containing the eigenvalues $\alpha_1,\ldots, \alpha_n$ and $\conj{\beta_1},\ldots, \conj{\beta_m}$ of $A$ and $B^H$.
    It holds:

    \begin{enumerate}[(i)]
        \item $\langle X, Y \rangle_{U,V}:= \langle U^{-1} X V^{-H}, U^{-1} Y V^{-H} \rangle_F$ is an inner product on $\mathbb{C}^{n\times m}$.
        \item ${(u_i {v_j}^H)}_{\substack{i=1,\ldots,n\\j=1,\ldots,m}}$ is an orthonormal basis of $\mathbb{C}^{n\times m}$ with respect to $\langle \cdot, \cdot \rangle_{U,V}$.
        \item The adjoint operator $\mathcal S^{\ast}: \Cmat{n}{m} \to \Cmat{n}{m}$ with respect to $\langle \cdot,\cdot \rangle_{U,V}$  \\
              is $\mathcal{S}^{\ast}(X)=U\overline{D_A} U^{-1}X  +  XV^{-H} D_{B^H}V^{H}$.
        \item $\mathcal S$ is a normal operator with respect to  $\langle \cdot, \cdot \rangle_{U,V}$.
    \end{enumerate}
\end{lemma}
\begin{proof}
    Note that $\langle \cdot, \cdot \rangle_{U,V}$ defines an inner product since
    \begin{align*}
        \langle u_i v_j^H , u_k v_l^H\rangle_{U,V} =
        \langle U^{-1} u_i  v_j^H V^{-H},  U^{-1} u_k  v_l^H V^{-H} \rangle_F =
        \langle e_i e_j^H , e_k e_l^H \rangle_F =
        \delta_{i,k}\delta_{j,l}.
    \end{align*}
    The matrices $u_i v_j^H \in \Cmat{n}{m}$ are orthogonal with respect to $\langle \cdot, \cdot \rangle_{U,V}$  and therefore linearly independent.
    Because  $\dim(\Cmat{n}{m})=n\cdot m$, the tuple ${(u_i v_j^H)}_{\substack{i=1,\ldots,n\\j=1,\ldots,m}}$ is an orthonormal basis of $\Cmat{n}{m}$. \\
    The following two computations show that $ \mathcal S$ commutes with its adjoint $\mathcal S^{\ast}$. \\
    Let $X,Y\in \Cmat{n}{m}$.
    \begin{align*}
        \langle S(X), Y \rangle_{U,V} & = \langle AX + XB , Y \rangle_{U,V}
        = \langle AX, Y \rangle_{U,V} + \langle XB , Y \rangle_{U,V}                                                                                      \\
                                      & = \langle U^{-1}AXV^{-H}, U^{-1}YV^{-H} \rangle_{F} + \langle U^{-1}XBV^{-H} , U^{-1}YV^{-H} \rangle_{F}          \\
                                      & = \langle D_A U^{-1}XV^{-H}, U^{-1}YV^{-H} \rangle_{F} +
        \langle U^{-1}XV^{-H}\overline{D_{B^H}} , U^{-1}YV^{-H} \rangle_{F}                                                                               \\
                                      & = \langle  U^{-1}XV^{-H}, \overline{D_A} U^{-1}YV^{-H}+ U^{-1}YV^{-H} D_{B^H} \rangle_{F}                         \\
                                      & = \langle  U^{-1}XV^{-H},U^{-1} \left(U\overline{D_A} U^{-1}Y  +  YV^{-H} D_{B^H}V^{H} \right) V^{-H} \rangle_{F} \\
                                      & = \langle  X,\left(U\overline{D_A} U^{-1}Y  +  YV^{-H} D_{B^H}V^{H} \right) \rangle_{U,V}                         \\
                                      & = \langle X, \mathcal S^{\ast}(Y)\rangle_{U,V}.
    \end{align*}
    Therefore, the adjoint of $\mathcal S$ is $\mathcal{ S}^{\ast}(X)=U\overline{D_A} U^{-1}X  +  XV^{-H} D_{B^H}V^{H}$.
    Moreover,
    \begin{align*}
        \mathcal S \mathcal S^{\ast}(X) & = \mathcal S (U\overline{D_A} U^{-1}X  +  XV^{-H} D_{B^H}V^{H})
        = \mathcal S (U\overline{D_A} U^{-1}X) + \mathcal S(XV^{-H} D_{B^H}V^{H})                                                       \\
                                        & = U D_A \overline{D_A}U^{-1} X +
        U\overline{D_A}U^{-1}XV^{-H}\overline{D_{B^H}}V^H                                                                               \\
                                        & \phantom{= } + U D_A U^{-1} X V^{-H} D_{B^H} V^{H} + X V^{-H}D_{B^H}\overline{D_{B^H}} V^{H}  \\
                                        & = U \overline{D_A} D_A U^{-1} X +
        U\overline{D_A}U^{-1}XV^{-H}\overline{D_{B^H}}V^H                                                                               \\
                                        & \phantom{= } + U D_A U^{-1} X V^{-H} D_{B^H} V^{H} + X V^{-H}\overline{D_{B^H}} D_{B^H} V^{H} \\
                                        & = \mathcal S^{\ast} (UD_A U^{-1}X) + \mathcal S^{\ast}(XV^{-H} \overline{D_{B^H}}V^{H})
        = \mathcal S^{\ast}\mathcal S(X).
    \end{align*}
    This means $\mathcal S$ and $\mathcal S^{\ast}$ commute and, therefore, by definition, $\mathcal S$ is normal.
\end{proof}

Now that we have an inner product on a Hilbert space for which $\mathcal S$ is normal, the second step is to compute the spectral decomposition of $\mathcal S$.
The spectral decomposition allows functional calculus and we get a formula for
$\mathcal S^{-1}$ and $e^{t\mathcal S}$. Since for normal operators, the operator norm is its spectral radius, we directly get a formula for the induced operator norm of $\mathcal S$.
We mention that in the case that $A$ and $B$ are unitarly diagonalizable, there is no need to exchange the
inner product as one can take the Frobenius inner product $\langle \cdot, \cdot \rangle_F$.

\begin{lemma}[Spectral Decomposition of $\mathcal S$]\label{satz:spectral_sylvester}
    Let the assumptions of \textup{Lemma~\ref{lemma:innerproduct}} hold.
    Then it holds:
    \begin{enumerate}[(i)]
        \item $\mathcal{S}(X)
                  = \summe{i=1}{n}\summe{j=1}{m}  (\alpha_i + \beta_{j}) \langle X, u_{i}v_{j}^H \rangle_{U,V} u_{i}v_{j}^H
                  = U\left( {\left( \alpha_i + \beta_j \right)}_{\substack{i=1,\ldots,n \\j =1,\ldots,m}}
                  \odot U^{-1} X V^{-H} \right) V^{H} $.
        \item $\|\mathcal S\|=\max\limits_{X\in \Cmat{n}{m}\setminus\{0\}} \frac{\|\mathcal S(X)\|_{U,V}}{\|X\|_{U,V}}=\max\limits_{i,j}|\alpha_i + \beta_j|$,
              where $\| X\|_{U,V} = \sqrt{\langle X,X\rangle_{U,V}}$.
        \item  $\mathcal{S}^{-1}(X) =U {\left( {\left(\frac{1}{\alpha_i + \beta_j}\right)}_{\substack{i=1,\ldots,n \\ j=1,\ldots,m}}\odot U^{-1}XV^{-H} \right)}V^{H}$
              and \newline
              $e^{t \mathcal S}(X) = U \left( {\left(e^{t(\alpha_i + \beta_j)}\right)}_{\substack{i=1,\ldots,n \\ j=1,\ldots,m}}\odot U^{-1}XV^{-H} \right)V^{H}$.
    \end{enumerate}
\end{lemma}
\begin{proof}
    \mbox{ }
    \begin{enumerate}[(i)]
        \item
              From  $AU=UD_A$ and $B^{H}V=VD_{B^{H}}$, we deduce $\mathcal{S}(u_k v_l^H) = A u_k v_l^H + u_k v_l^H B = (\alpha_k + \beta_l) u_k v_l^H$.
              Representing $\mathcal{S}(X)\in \Cmat{n}{m}$ as well as $X\in \Cmat{n}{m}$ as a Fourier series and exploiting linearity of $\mathcal S$ and
              $\langle \cdot, \cdot \rangle_{U,V}$ yields
              \begin{align*}
                  \mathcal{S}(X) & = \summe{i=1}{n}   \summe{j=1}{m}
                  \langle  \mathcal{S}(X),u_i v_j^H \rangle_{U,V}  u_i v_j^H
                  = \summe{i=1}{n}   \summe{j=1}{m}
                  \langle \mathcal{S}(\summe{k=1}{n}\summe{l=1}{m} \langle  X ,u_k v_l^H\rangle_{U,V} u_k v_l^H),
                  u_i u_j^H
                  \rangle_{U,V}u_i v_j^H                                                                \\
                                 & = \summe{i,k=1}{n} \summe{j,l=1}{m}
                  \langle X, u_k v_l^H  \rangle_{U,V}
                  \langle \mathcal{S}(u_k v_l^H), u_i v_j^H \rangle_{U,V} u_i v_j^H                     \\
                                 & = \summe{i,k=1}{n} \summe{j,l=1}{m}
                  (\alpha_k + \beta_l)
                  \langle X,u_k v_l^H \rangle_{U,V}
                  \langle u_k v_l^H ,u_i v_j^H\rangle_{U,V} u_i v_j^H                                   \\
                                 & = \summe{i=1}{n} \summe{j=1}{m}
                  (\alpha_i + \beta_j)
                  \langle X,u_i v_j^H \rangle_{U,V} u_i v_j^H
                  = U
                      {\left((\alpha_i + \beta_j)
                          \langle X,u_i v_j^H \rangle_{U,V}
                          \right)}_{\substack{i=1,\ldots,n                                              \\ j=1,\ldots,m}}
                  V^{H}                                                                                 \\
                                 & = U
                  \left(
                  {\left(\alpha_i + \beta_j\right)}_{\substack{i=1,\ldots,n                             \\ j=1,\ldots,m}}
                  \odot
                  {\left( \langle  U^{-1}XV^{-H},e_i e_j^H \rangle_{F} \right)}_{\substack{i=1,\ldots,n \\ j=1,\ldots,m}}
                  \right)
                  V^{H}                                                                                 \\
                                 & = U
                  \left(
                  {\left(\alpha_i + \beta_j\right)}_{\substack{i=1,\ldots,n                             \\ j=1,\ldots,m}}
                  \odot
                  U^{-1}XV^{-H}
                  \right)
                  V^{H}.
              \end{align*}
        \item The claim about the norm follows from a direct application of the fundamental functional analytical result on compact normal operators, see, e.g.,~\cite[Thm.~\rz{VI}.3.2]{Wer00}.%
        \item With the spectral decomposition of $\mathcal S$ one can resort to functional calculus, cf.~\cite[Cor.~9.3.38]{Pal01},~\cite[Kor.~\rz{IX}.3.8]{Wer00}, and obtain the formula for $\mathcal S^{-1}$ under the additional assumption that $\alpha_i + \beta_j \neq 0$.
    \end{enumerate}
\end{proof}

Using the spectral decomposition and functional calculus we find that, under the assumptions of Lemma~\ref{lemma:innerproduct}, the solution of the differential Sylvester equation
\begin{align*}
    \dot{X}(t) & =  AX(t) +X(t)B + C = \mathcal{S}(X(t)) +  C, \\
    X(0)       & = D,
\end{align*}
has the form
\begin{align}
    X(t) & =  e^{t\mathcal S}(D) + \int\limits_{0}^{t}e^{(t-s)\mathcal S}(C) \dx{s}             \notag \\
         & =  U \left( {\left(e^{t(\alpha_i + \beta_j)}\right)}_{\substack{i=1,\ldots,n                \\ j=1,\ldots,m}}\odot U^{-1}DV^{-H} \right)V^{H}    +
    \int\limits_{0}^{t}
    U \left( {\left(e^{(t-s)(\alpha_i + \beta_j)}\right)}_{\substack{i=1,\ldots,n                      \\ j=1,\ldots,m}}\odot U^{-1}CV^{-H} \right)V^{H} \dx{s}         \notag \\
         & =  U \left(
    {\left(e^{t(\alpha_i + \beta_j)}\right)}_{\substack{i=1,\ldots,n                                   \\ j=1,\ldots,m}} \odot U^{-1}DV^{-H}
    +
    {\left( \int\limits_{0}^{t} e^{(t-s)(\alpha_i + \beta_j)} \dx{s}\right)}_{\substack{i=1,\ldots,n   \\ j=1,\ldots,m}}\odot U^{-1}CV^{-H}
    \right)
    V^{H}, \label{eq:spectrldecsolformula}
\end{align}
with the involved scalar integrals given explicitly as:
\begin{align*}
    \integraltwo{0}{t} e^{(t-s){(\alpha_i +\beta_j)}}\dx{s} =
    \begin{cases}
        \frac{e^{t(\alpha_i + \beta_j)}-1}{\alpha_i + \beta_j} & \mbox{if } \alpha_i + \beta_j\neq 0, \\
        t                                                      & \mbox{if } \alpha_i + \beta_j= 0
    \end{cases}.
\end{align*}

\section{Variation of Constants}\label{sec_variation_const}
The application of the variation of constants formula leads to yet another solution formula for the SLE~\eqref{eq:introsylvestereq}.
\begin{lemma}[Variations of Constants, {\cite[Ch.~13]{Bro70}}]\label{lemma:variations_of_constants}
    Let $A\in \Cmat{n}{n}, B \in \Cmat{m}{m}, C\in \Cmat{m}{n}, D\in \Cmat{m}{n}$ with $\spektrum{A}\cap \spektrum{-B}=\emptyset$.
    The differential Sylvester equation
    \begin{align*}
        \dot{X}(t) & = AX(t) + X(t)B + C=  \mathcal{S}(X(t))  + C, \\
        X(0)       & =D,
    \end{align*}
    has the solution
    \begin{equation}\label{eq:variatconstSylvSplit}
        X(t) = e^{t\mathcal S}(D) + \mathcal{S}^{-1}(-C)   - e^{t\mathcal S}\mathcal{S}^{-1}(-C).
    \end{equation}
\end{lemma}

\begin{proof}
    Because of $\spektrum{A}\cap \spektrum{-B}= \emptyset$, the inverse $\mathcal{S}^{-1}$ exists and we can rewrite the solution formula~\eqref{eq:variatconstSylv} as
    \begin{align*}
        X(t) & =e^{t\mathcal S}(D) + \integraltwo{0}{t}e^{(t-s)\mathcal S}(C) \dx{s}
        =e^{t\mathcal S}(D) +  {\left[-\mathcal S^{-1}e^{(t-s)\mathcal S}(C)\right]}_{s=0}^{s=t} \\
             & =e^{t\mathcal S}(D) + \mathcal S^{-1}(-C) - \mathcal S^{-1}e^{t\mathcal S}(-C)
    \end{align*}
    and confirm that
    $ X(0) = D  + \mathcal{S}^{-1}(-C) - \mathcal S^{-1}(-C)= D$.
\end{proof}

From formula~\eqref{eq:variatconstSylvSplit}, we find that the solution can be written as the solution of the algebraic Sylvester equation and a time dependent part. We will make use of this fact in the numerical scheme that we propose in Section~\ref{sec_feasible_numerical_appr}.
%

\section{Solution as Taylor Series}\label{sec_taylor}
In this section we use Taylor series to derive a solution formula. From this we can read off suitable Krylov subspaces
for our projection approach in the next section.
\begin{lemma}[Taylor Series Solution Formula]
    Let $A\in \Cmat{n}{n}, B \in \Cmat{m}{m}, C\in \Cmat{m}{n}, D\in \Cmat{m}{n}$.
    The differential Sylvester equation
    \begin{align*}
        \dot{X}(t) & = AX(t) + X(t)B + C =\mathcal{S}(X(t))  + C,  \label{eq:sylvester_autonom} \\
        X(0)       & =D,
    \end{align*}
    has the unique solution
    \begin{equation}\label{eq:solastayl}
        X(t) = D + \summe{k=1}{\infty}\frac{t^k}{k!}(\mathcal S^k(D)+\mathcal S^{k-1}(C)).
    \end{equation}
\end{lemma}
\begin{proof}
    \begin{align*}
        \norm{D}                    +   \summe{k=1}{\infty}\norm{\frac{t^k}{k!}(\mathcal S^k(D)+\mathcal S^{k-1}(C))}
         & \leq     \norm{D}                    +   \summe{k=1}{\infty}\frac{|t|^k}{k!}(\norm{\mathcal S^k(D)}+\norm{\mathcal S^{k-1}(C)})              \\
         & \leq     \norm{D}                    +   \summe{k=1}{\infty}\frac{|t|^k}{k!}( \norm{\mathcal S}^k\norm{D} + \norm{\mathcal S}^{k-1}\norm{C}) \\
         & =    \norm{D}e^{|t|\norm{\mathcal{S}}}   +   \norm{C}\summe{k=0}{\infty}\frac{|t|^{k+1}}{(k+1)!} \norm{\mathcal S}^{k}                       \\
         & \leq     \norm{D}e^{|t|\norm{\mathcal{S}}}   +   |t|\norm{C}e^{|t|\norm{\mathcal{S}}}.
    \end{align*}
    The series converges absolutely and since $(\mathbb{C}^{n\times m},\norm{\cdot})$ is a Banach space, the series converges for every $t \in \Rfield$. Therefore its radius of convergence is infinity,
    $X$ is continuously differentiable and can be differentiated term-wise.
    Since, furthermore,
    \begin{align*}
        X(0)      & =D  \intertext{and}
        \dot X(t) & =
        \summe{k=1}{\infty}\frac{t^{k-1}}{(k-1)!}(\mathcal{S}^{k}(D)+ \mathcal{S}^{k-1}(C))               \\ &=
        \summe{k=0}{\infty}\frac{t^{k}}{k!}(\mathcal{S}^{k+1}(D)+ \mathcal{S}^{k}(C)) =
        \mathcal{S}(D) + \summe{k=1}{\infty}\frac{t^{k}}{k!}(\mathcal{S}^{k+1}(D)+ \mathcal{S}^{k}(C)) +C \\ &=
        \mathcal{S}(D  + \summe{k=1}{\infty}\frac{t^{k}}{k!}(\mathcal{S}^{k}(D)+ \mathcal{S}^{k-1}(C))) +C =
        \mathcal{S}(X(t)) + C,
    \end{align*}
    $X(t)$ is the unique solution.
\end{proof}

If we assume that $D$ and $C$ are given in factored form, then we can exploit this to rewrite the truncated series in a closed form of a matrix product. 

\begin{remark}\label{remark:taylor}
    Let $D=D_1D_2^H$ and $C=C_1C_2^H$ with
    $D_1 \in \Cmat{n}{d}, D_2 \in \Cmat{m}{d},C_1 \in \Cmat{n}{c}$ and $C_2 \in \Cmat{m}{c}$.
    Then, having truncated the two parts of the series~\eqref{eq:solastayl} after $m_1$ and $m_2$ summands, respectively, we can rewrite the solution approximation as
    \begin{align*}
        X_{m_1,m_2}(t) & =  \summe{k=0}{m_1}    \frac{t^k}{k!}\mathcal S ^k (D)             + \summe{k=1}{m_2} \frac{t^k}{k!} \mathcal S^{k-1}(C)                                     \\
                       & =  \summe{k=0}{m_1}    \frac{t^k}{k!}\mathcal S ^k (D_1D_2^H)          + \summe{k=1}{m_2} \frac{t^k}{k!} \mathcal S^{k-1}(C_1C_2^H)                          \\
                       & =  \summe{k=0}{m_1}    \frac{t^k}{k!}{(\mathcal H + \mathcal V)}^k (D_1D_2^H)  + \summe{k=1}{m_2} \frac{t^k}{k!} {(\mathcal H + \mathcal V)}^{k-1}(C_1C_2^H) \\
                       & =  \summe{k=0}{m_1}\summe{i=0}{k}\frac{t^k}{k!}\binom{k}{i}\mathcal H^{k-i}\mathcal V ^{i}(D_1D_2^H)   +
        \summe{k=1}{m_2}\summe{i=0}{k-1}\frac{t^k}{k!}\binom{k-1}{i}\mathcal H^{k-1-i}\mathcal V ^{i}(C_1 C_2^H)                                                                      \\
                       & =
        \summe{k=0}{m_1}\summe{i=0}{k}  \frac{t^k}{k!}\binom{k}{i} A^{k-i}D_1D_2^{H}B^i +
        \summe{k=1}{m_2}\summe{i=0}{k-1}\frac{t^k}{k!}\binom{k-1}{i} A^{k-1-i}C_1 C_2^{H}B^i,
    \end{align*}
    with the explicit representation of the sums
    \begin{align*}
         & \summe{k=0}{m_1}\summe{i=0}{k}\frac{t^k}{k!}\binom{k}{i}A^{k-i} D_1 D_2^H B^{i}
        =                                                                                  \\
         & \begin{bmatrix} D_1, A D_1, \ldots, A^{m_1}D_1 \end{bmatrix}
        \left(
        \begin{bmatrix}     \frac{t^0}{0!}\binom{0}{0}         & \frac{t^1}{1!}\binom{1}{1} & \frac{t^2}{2!}\binom{2}{2} & \cdots  & \frac{t^m_1}{m_1!}\binom{m_1}{m_1} \\
            \frac{t^1}{1!}\binom{1}{0}         & \iddots                    & \iddots                    & \iddots & 0                                  \\
            \frac{t^2}{2!}\binom{2}{0}         & \iddots                    & \iddots                    & \iddots & \vdots                             \\
            \vdots                             & \iddots                    & \iddots                    & \iddots & 0                                  \\
            \frac{t^{m_1}}{m_1!}\binom{m_1}{0} & 0                          & \cdots                     & 0       & 0
        \end{bmatrix}
        \kron E_{d,d}\right)
        \begin{bmatrix}  D_2^H \\ D_2^H B  \\ \vdots \\ D_2^H B^{m_1} \end{bmatrix}
    \end{align*}
    and
    \begin{align*}
         & \summe{k=1}{m_2}\summe{i=0}{k-1} \frac{t^k}{k!}\binom{k-1}{i} A^{k-1-i}C_1 C_2^{H}B^i = \\
         & \begin{bmatrix} C_1, A C_1,\ldots,A^{m_2-1}C_1  \end{bmatrix}
        \left(
        \begin{bmatrix}
            \frac{t^1}{1!}\binom{0}{0}             & \frac{t^2}{2!}\binom{1}{1} & \frac{t^3}{3!}\binom{2}{2} & \cdots  & \frac{t^{m_2}}{m_2!}\binom{m_2-1}{m_2-1} \\
            \frac{t^2}{2!}\binom{1}{0}             & \iddots                    & \iddots                    & \iddots & 0                                        \\
            \frac{t^3}{3!}\binom{2}{0}             & \iddots                    & \iddots                    & \iddots & \vdots                                   \\
            \vdots                                 & \iddots                    & \iddots                    & \iddots & 0                                        \\
            \frac{t^{m_2}}{m_2!}\binom{{m_2}-1}{0} & 0                          & \cdots                     & 0       & 0
        \end{bmatrix}
        \kron E_{c,c}\right)
        \begin{bmatrix}  C_2^H \\ C_2^H B  \\ \vdots  \\ C_2^H B^{m_2-1} \end{bmatrix}.
    \end{align*}
\end{remark}

\section{Feasible Numerical Solution Approaches}\label{sec_feasible_numerical_appr}
In this section, we briefly note that, for various reasons, all presented solution representations are not feasible for a straight-forward numerical approximation, in particular in a large-scale sparse setting.

A common reason is that none of the formulas supports a sparse representation of the solutions such that exorbitant amounts of memory will be required.

Limitations in memory will doubly affect the solution representation through the spectral decomposition~\eqref{eq:spectrldecsolformula} since also the basis matrices $U$ and $V$ are generically dense matrices. Apart from that, the computation of a spectral decomposition is typically computationally expensive and can be ill conditioned. Nonetheless, the spectral decomposition formula is useful to construct exact solutions for given coefficients with known spectral decompositions.

Another issue is the unfeasible computation of the full matrix exponential in all variants~\eqref{eq:varofconstsAXXB},~\eqref{eq:variatconstSylv}, and~\eqref{eq:variatconstSylvSplit} of the \emph{variation of constants} formula.
A possible remedy is the approximation the action of the matrix exponential on a low-rank matrix in a Krylov subspace.

The approach to the solution via a Taylor series (see Section~\ref{sec_taylor}) seems best suited for the large-scale case since, at least in the symmetric case, the formulas provided in Remark~\ref{remark:taylor} allow for a solution representation in factored form with the original coefficients. One problem
here is that the truncated Taylor series only leads to good approximations locally around the point of expansion.

We will, however, exploit and combine certain parts of the solution representations to propose an algorithm for fast and memory efficient solution approximations.

We consider the stable, linear time-invariant case, i.e.~we assume that $A\in \Rmat{n}{n}$, $B\in \Rmat{n}{p}$, and $\spektrum{A}\subseteq \C^{-}$.
We consider the differential Lyapunov equation
\begin{subequations}\label{eqn:autodle}
    \begin{align}
        \dot{X}(t) & = A^T X(t) + X(t)A  + BB^T,    \label{eqn:autodle_eqn} \\
        X(0)       & = 0.                           \label{eqn:autodle_ini}
    \end{align}
\end{subequations}

By Lemma~\ref{lemma:variations_of_constants}, we have that the solution splits into a constant part and a time dependent part. From Remark~\ref{remark:taylor}, we infer that the solution is contained in a Krylov subspace. We combine both observations in the following numerical solution approach:

\textbf{1. Factors of the time constant Part as Krylov space basis}
With $A$ stable, the associated algebraic Lyapunov $A^T X + X A + BB^T=0$ has a unique symmetric positive semi-definite solution $X_{\infty}$ that can be written in factored form $X_\infty = Z_{\infty}Z_{\infty}^T$, with $Z_{\infty} \in \Rmat{n}{q}$ and $\rank(Z_{\infty})=q\leq n$. Moreover, since $A$ is stable, it holds (see~\cite[Ch.~13]{Bro70}) that
\begin{align*}
    \range(Z_{\infty})=\range(X_{\infty})=\range\left( \left[B,A^{T}B, \ldots, {(A^T)}^{n-1} B \right] \right).
\end{align*}

\textbf{2. Factors of the time dependent part evolve in the same Krylov space}
With $X(t)= X_{\infty} + \tilde{X}(t)$, we obtain that
\begin{align*}
    \dot{\tilde{X}}(t) & = A^T \tilde{X}(t)  + \tilde{X} (t)A, \\
    \tilde{X}(0)       & = -X_{\infty},
\end{align*}
where $\tilde{X}(t)$ is given by $\tilde{X}(t)= -e^{tA^T}X_{\infty} e^{tA} =-e^{tA^T}Z_{\infty} Z_{\infty}^T e^{tA}=:-\tilde{Z}(t)\tilde{Z}{(t)}^T$.

\textbf{3. Orthogonalize the basis and compute the time dependent factors}
By means of the singular value decomposition of $Z_{\infty}$, we obtain an orthogonal matrices $Q_{\infty}$ and $V_{\infty}$
with $\range(Q_{\infty}) = \range(Z_{\infty})$,
with $Z_{\infty} = Q_{\infty}S_{\infty} V_{\infty}^T$, and with $Q_{\infty}\in \Rmat{n}{q},\ S_{\infty}\in \Rmat{q}{q},\ V_{\infty}\in \Rmat{q}{q}$.
Like $Z_\infty$, the columns of $Q_{\infty}$ span an $A^T$ invariant subspace and with $Q_{\infty}^T Q_{\infty}=E_{q,q}$ it holds that
\begin{align*}
    A^T Q_{\infty}      & = Q_{\infty} Q_{\infty}^T A^T Q_{\infty},       \\
    e^{tA^T} Q_{\infty} & = Q_{\infty}  e^{tQ_{\infty}^T A^T Q_{\infty}},
\end{align*}
from which we confer that
\begin{align*}
    \tilde{Z}(t){\tilde{Z}(t)}^T & = \left(e^{tA^T}Z_{\infty}\right) {\left(e^{tA^T}Z_{\infty}\right)}^T = \left(e^{tA^T}Q_{\infty} S_{\infty}\right) {\left(e^{tA^T}Q_{\infty} S_{\infty}\right)}^T \\
                                 & = \left(Q_{\infty}e^{tQ_{\infty}^{T}A^{T}Q_{\infty}} S_{\infty}\right){\left(Q_{\infty}e^{tQ_{\infty}^{T}A^{T}Q_{\infty}} S_{\infty}\right)}^T.
\end{align*}
We define ${z}(t)=e^{tQ_{\infty}^{T}A^{T}Q_{\infty}} S_{\infty}\in \Rmat{q}{q}$ and find that ${z}$ can be obtained by solving
\begin{subequations}\label{eq:projected_td_factors}
    \begin{align}
        \dot{{z}}(t) & = Q_{\infty}^{T}A^{T}Q_{\infty} {z}(t) \\
        {z}(0)       & = S_{\infty},
    \end{align}
\end{subequations}
which is a matrix valued ODE that can be solved column-wise or by computing the matrix exponential $e^{tQ_{\infty}^{T}A^{T}Q_{\infty}}$.

The  solution of the differential Lyapunov equation is, thus, given by
$X(t) = Z_{\infty} Z_{\infty}^T - Q_{\infty}{z}(t){z}{(t)}^T Q_{\infty}^T$. \\

\begin{remark}
    The differential equation for $z$ is of size $q\times q$, which can be much smaller than $n$, if the solution of the algebraic Lyapunov equation $X_{\infty}=Z_{\infty}Z_{\infty}^T$ has low-rank. Moreover, the orthogonalization of the basis allows for the detection of the numerical rank and for a compression of $Z_\infty$ through truncating singular values that are smaller than a certain threshold.

    We further note that, with minor adjustments, all arguments also hold for the generalized differential Lyapunov equation
    \begin{align*}
        M^T\dot{X}(t)M & = A^T X(t) M + {M}^T X(t)A  + BB^T, \\
        X(0)           & = 0,
    \end{align*}
    with $M\in \Rmat{n}{n}$ nonsingular
    that can accommodate, e.g., a mass matrix from a finite element discretization.
\end{remark}

In summary, the proposed approach reads as written down in Algorithm~\ref{alg:projection}.
\begin{algorithm}
    \caption{Projection approach for generalized differential Lyapunov equations}\label{alg:projection}
    \begin{algorithmic}[1]                  
        \REQUIRE{} $M,A\in \Rmat{n}{n}$ with $\spektrum{AM^{-1}}\subseteq \Cfield_{-}$ and $B\in \Rmat{n}{p}$.
        \ENSURE{} $X(t) = Z_{\infty}Z_{\infty}^T - Q_{\infty}{z}(t){z}{(t)}^T Q_{\infty}^T$
        that approximates the solution to \\ $M^T\dot{X}(t)M=A^T X(t) M + M^{T}X(t) A + BB^T,\ X(0)=0$.
        \vspace{.1in}
        \STATE{} \% Solve Lyapunov equation:
        \STATE{} $A^T X_{\infty}M + M^{T}X_{\infty} A = -BB^T$ for $X_{\infty} \approx Z_{\infty} Z_{\infty}^T$ and $Z_{\infty}\in \Rmat{n}{q}$.
        \vspace{.1in}
        \STATE{} \% Compute singular value decomposition:
        \STATE{} $[Q_{\infty}, S_{\infty},\sim] = \text{svd}(Z_{\infty},0)$.
        \vspace{.1in}
        \STATE{} \% Set tolerance to largest singular value times machine epsilon:
        \STATE{} $tol = \vareps_{\text{machine}} \cdot S_{\infty}(1,1)$.
        \vspace{.1in}
        \STATE{} \% Truncate all singular values smaller than tolerance and get truncated low-rank factor:
        \STATE{} $idx = \text{diag}(S_{\infty})\geq tol$.
        \STATE{} $S_{\infty} \leftarrow S_{\infty}(idx,idx)$.
        \STATE{} $Q_{\infty} \leftarrow Q_{\infty}(:,idx)$.
        \STATE{} $Z_{\infty} \leftarrow Q_{\infty}S_{\infty}$.
        \vspace{.1in}
        \STATE{} \% Compute projected system and solve:
        \IF{$M$ is symmetric positive definite}
        \STATE{} $M_F = Q_{\infty}^T M^T Q_{\infty}$.
        \STATE{} $A_F = Q_{\infty}^T A^T Q_{\infty}$.
        \ELSE{}
        \STATE{} $M_F = E$.
        \STATE{} $A_F = Q_{\infty}^T M^{-T}A^T Q_{\infty}$.
        \ENDIF{}
        \vspace{.1in}
        \FOR{k=1,\ldots,\text{ cols ($S_{\infty}$)}}
        \STATE{} Solve: $M_F\dot{{z}}(:,k)(t)=A_F{z}(:,k)(t),\ {z}(:,k)(0)={S_{\infty}}(:,k)$.
        \ENDFOR{}
    \end{algorithmic}
\end{algorithm}
\pagebreak

\section{Numerical Results}\label{sec_numerical_res}
\subsection{Setup}\label{sec_numerical_res_setup}
To quantify and illustrate the performance of Algorithm~\ref{alg:projection}, we consider differential Lyapunov equations that are used to define optimal controls for a finite element discretization of a heat equation; see~\cite{morBenS05} for the model description.
Namely, we solve the differential Lyapunov equations:
\begin{equation}\label{eq:numtestDLEc}
    M \dot{X}(t) M^T = A X(t) M^T + M X(t) A^T + BB^T,      \quad X(0)  = 0. \tag{DLE--1}
\end{equation}
\begin{equation}\label{eq:numtestDLEo}
    M^T \dot{X}(t) M  = A^T X(t) M + M^T X(t) A + C^T C,    \quad X(0)  = 0. \tag{DLE--2}
\end{equation}
that are defined through matrices $M$, $A\in \Rmat{n}{n}$ that are symmetric, $M$ is positive definite and $A$ stable, $B\in \Rmat{n}{7}$ and $C\in \Rmat{6}{n}$.
For computing the error, we precomputed the spectral decomposition of $(A,M)$ and constructed the exact solution by means of the formula from Lemma~\ref{satz:spectral_sylvester}.
The memory consuming computation of the spectral decomposition was done on a compute server with  4 $\times$ \xeon{} CPU E7--8837 @ 2.67GHz with 8 cores and 1 TB Ram and \matlab{} 2015b.
All other computations were carried out on a machine with 2 $\times$ \xeon{} CPU E5--2640 v3 @ 2.60GHz with 8 Cores and 64 GB Ram and \matlab{} 2017a. We have used the low-rank ADI iteration
implemented in
\mexmess{}\cite{messweb}~to solve the algebraic Lyapunov equations; as required for Algorithm~\ref{alg:projection} (Step 2). \\
We solve the resulting projected ODE system column-wise using~\matlab{} ODE solvers
\texttt{ode45}, \texttt{ode23}, \texttt{ode113}, \texttt{ode15s}, \texttt{ode23s}, \texttt{ode23t} and \texttt{ode23tb}, with the parameters for the $\texttt{odeset}$ function set as follows:
\begin{multicols}{3}
    \begin{itemize}
        \item \texttt{RelTol}: $1\cdot 10^{-9}$
        \item \texttt{AbsTol}:  $1\cdot 10^{-10}$
        \item \texttt{Stats}: \texttt{off}
        \item \texttt{NormControl}: \texttt{off}
        \item \texttt{BDF}: \texttt{on}
        \item \texttt{Jacobian}: $A_F$
        \item \texttt{JPattern}: \texttt{logical} ($A_F$)
        \item \texttt{Mass}: $M_F$
        \item \texttt{MassSingular}: \texttt{no}
        \item \texttt{MStateDependence}: \texttt{none}
    \end{itemize}
\end{multicols}
As the time interval, we considered $[0,4500]$.

\subsection{Projection Approach}\label{sec_numerical_res_proj}
The initial step of Algorithm~\ref{alg:projection} requires the solutions to the associated algebraic Lyapunov equations. For this task we call \mexmess~that iteratively computes the solutions up to the following absolute and relative residuals
\\
\begin{align*}
    \norm{A Z_{\infty} Z_{\infty}^T M^T + M Z_{\infty} Z_{\infty}^T A^T + BB^T}_2\ \text{or}\ \norm{ A^T Z_{\infty} Z_{\infty}^T M + M^T Z_{\infty} Z_{\infty}^T A +  C^{T}C}_2.
\end{align*}
and
\begin{align*}
    \frac{\norm{A Z_{\infty} Z_{\infty}^T M^T + M Z_{\infty} Z_{\infty}^T A^T + BB^T}_2}{\norm{BB^T}_2}\ \text{or}\ \frac{\norm{A^T Z_{\infty} Z_{\infty}^T M +  M^T  Z_{\infty} Z_{\infty}^T A + C^{T}C}_2}{\norm{C^{T}C}_2}.
\end{align*}
The achieved values for the different test setups as well as the number of columns of the corresponding $Z_{\infty}$ after truncation (see Step 7 of Algorithm~\ref{alg:projection}), that define the dimension of the reduced model for the time dependent part, are listed in Tables~\ref{tab_res_bbt} and~\ref{tab_res_ctc}.

\begin{table}[htb]
    \centering
    \begin{tabular}{|c|c|c|c|}
        \hline
        size  & size of $z(t)$   & absolute residual         & relative residual         \\ \hline
        1357  & $261 \times 261$ & $3.413488 \cdot 10^{-25}$ & $7.748357 \cdot 10^{-12}$ \\ \hline
        5177  & $302 \times 302$ & $1.037846 \cdot 10^{-25}$ & $4.728703 \cdot 10^{-12}$ \\ \hline
        20209 & $376 \times 376$ & $6.053185 \cdot 10^{-26}$ & $5.525974 \cdot 10^{-12}$ \\ \hline
    \end{tabular}
    \caption{Residuals for $A X M^T + M X A^T + BB^T = 0$.}\label{tab_res_bbt}
\end{table}

\pagebreak
\begin{table}[htb]
    \centering
    \begin{tabular}{|c|c|c|c|}
        \hline
        size  & size of $z(t)$   & absolute residual         & relative residual         \\ \hline
        1357  & $230 \times 230$ & $1.011843 \cdot 10^{-10}$ & $8.432027 \cdot 10^{-12}$ \\ \hline
        5177  & $259 \times 259$ & $5.595100 \cdot 10^{-11}$ & $4.662583 \cdot 10^{-12}$ \\ \hline
        20209 & $310 \times 310$ & $4.382439 \cdot 10^{-11}$ & $4.382439 \cdot 10^{-12}$ \\ \hline
    \end{tabular}
    \caption{Residuals for $A^T X M + M^T X A +  {C}^T C = 0$.}\label{tab_res_ctc}
\end{table}

We report the absolute and relative errors
\begin{align*}
    \norm{X(t) - X_{ref}(t)}_2\quad \text{and}\quad \frac{\norm{X(t) - X_{ref}(t)}_2}{\norm{X_{ref}(t)}_2},
\end{align*}
where $X$ is the numerical solution obtained from Algorithm~\ref{alg:projection} with various ODE solvers and where the reference solution $X_{ref}$ was obtained from spectral decomposition of $(A,M)$.

We plot the numerical errors and $\norm{X_{ref}(t)}_2$ on the initial short time interval $[0,10]$, where most of the evolution is happening, and on the full time interval $[0,4500]$ in Appendix~\ref{sec_appendix_numerical_res_proj},
Figures
\ref{fig:appendix:proj:rel:long:1357:none}--\ref{fig:appendix:proj:abs:short:1357:none},
\ref{fig:appendix:proj:rel:long:5177:none}--\ref{fig:appendix:proj:abs:short:5177:none},
\ref{fig:appendix:proj:rel:long:20209:none}--\ref{fig:appendix:proj:abs:short:20209:none},
\ref{fig:appendix:proj:rel:long:1357:trans}--\ref{fig:appendix:proj:abs:short:1357:trans},
\ref{fig:appendix:proj:rel:long:5177:trans}--\ref{fig:appendix:proj:abs:short:5177:trans} and
\ref{fig:appendix:proj:rel:long:20209:trans}--\ref{fig:appendix:proj:abs:short:20209:trans}.

In view of the performance of the different ODE solvers, we can interprete the presented numbers and plots as follows: the solver $\texttt{ode15s}$, which is a stiff solver of variable order, performs best in time and accuracy. Due to the stiffness, the error is oscillating for the solvers \texttt{ode45}, \texttt{ode23}, \texttt{ode113}. Note that the discrete Laplacian that is encoded in the coefficient matrix $A$ becomes stiffer with a finer space discretization, i.e.~for larger $n$. Accordingly, the computational times for the non-stiff solvers grow with $n$ at a higher rate than the stiff solvers; see Figure~\ref{fig:odesolvtimings}.

The solution of~\eqref{eq:numtestDLEo} itself is large in norm what makes the relative error stagnate around the prescribed tolerance and the absolute error comparatively large; see the plots in Appendix~\ref{sec_appendix_numerical_res_proj_dleo}. Particularly, the non-stiff solvers (and \texttt{ode15s}) achieve this error level and the oscillations due to stiffness are dominated by the approximation error. This might be the reason, why for~\eqref{eq:numtestDLEo}, despite the fact that the coefficient matrices are still stiff, the non-stiff solvers perform better than the stiff solvers (except \texttt{ode15s}). Nonetheless, again the computation times for the non-stiff solver grow at a higher rate with the increasing stiffness that comes with increasing $n$; see Figure~\ref{fig:odesolvtimings}.

Because $X(t) \rightarrow 0$ for $t \searrow 0$,
the plots for the relative error spread out for small times.

Except \texttt{ode15s}, there is no general rule which solver performs better in terms of computational time; see Figure~\ref{fig:odesolvtimings}. The timings may change, when different relative and absolute error tolerances are used in the \matlab{} \texttt{odeset} function.

Finally, we want to make the following remark: instead of integrating the projected ODE~\eqref{eq:projected_td_factors} which is linear with constant coefficients of moderate dimension, one may consider using the \emph{Schur-Parlett}~\cite[Ch. 10]{Hig08} algorithm to compute the matrix exponential.
The initialization efforts of the \emph{Schur-Parlett} algorithm will pay off, if the matrix exponential has to be evaluated for many different values of $t$. Also, asymptotically, the storage requirements for $z(t)$ will be lifted, since the matrix exponential for a given $t$ can be computed on demand.
Nonetheless, we used the ODE approach, which we think is more efficient because of the sophisticated \matlab{} ODE solver implementations that come, e.g, with step size selection methods integrated.

The code of the implementation with the precomputed spectral decompositions needed to construct the exact solution is available as mentioned in Figure~\ref{fig:linkcodndat}.
\begin{figure}[h!]
    \begin{framed}
        \textbf{Code and Data Availability} \\
        The source code of the implementations used to compute the presented results is available from:
        \begin{center}
            \href{https://doi.org/10.5281/zenodo.1484327}{\texttt{doi:10.5281/zenodo.1484327}}
            \href{https://gitlab.mpi-magdeburg.mpg.de/behr/diff_lyap_eqn_solution_formulas_release}{\texttt{https://gitlab.mpi-magdeburg.mpg.de/behr/diff\_lyap\_eqn\_solution\_formulas\_release}}
        \end{center}
        under the GPLv2+ license
        and is authored by Maximilian Behr.
    \end{framed}
    \caption{Link to code and data.}\label{fig:linkcodndat}
\end{figure}

\pagebreak
\subsection{Computational Time}\label{sec_numerical_res_comptime}

\pgfplotsset{select coords between index/.style 2 args={
            x filter/.code={
                    \ifnum\coordindex<#1\def\pgfmathresult{}\fi
                    \ifnum\coordindex>#2\def\pgfmathresult{}\fi
                }
        }
}

\pgfplotsset{
    discard if not/.style 2 args={
            x filter/.append code={
                    \edef\tempa{\thisrow{#1}}
                    \edef\tempb{#2}
                    \ifx\tempa\tempb
                    \else
                        \def\pgfmathresult{inf}
                    \fi
                }
        }
}

\pgfplotsset{
    walltimestyle_proj/.style  = {
            yscale=1,
            xscale=0.75,
            ymin=1e0,
            ymax=5e4,
            extra y ticks={5, 5e1, 5e2, 5e3, 5e4},
            bar width=2mm,
            ylabel={Computational Time in Seconds},
            ylabel style={align=center,yshift=6mm},
            xlabel style={align=center,yshift=0mm},
            symbolic x coords={1357, 5177, 20209},
            enlargelimits=0.25,
            enlarge y limits=auto,
            xtick=data,
            nodes near coords align={horizontal},
        },
}

\pgfplotsset{
    walltimestyle_bdf/.style  = {
            yscale=1,
            xscale=0.75,
            ymin=1e0,
            ymax=5e4,
            extra y ticks={5, 5e1, 5e2, 5e3, 5e4},
            bar width=2mm,
            ylabel={Computational Time in Seconds},
            ylabel style={align=center,yshift=6mm},
            xlabel style={align=center,yshift=0mm},
            symbolic x coords={{k4}, {k6}, {k8}},
            xtick = {{k4}, {k6}, {k8}},
            xticklabels={$2^{-4}$, $2^{-6}$, $2^{-8}$},
            enlargelimits=0.25,
            enlarge y limits=auto,
            nodes near coords align={horizontal},
        },
}

\begin{minipage}{0.5\textwidth}
    \centering
    \begin{tikzpicture}[baseline]
        \begin{axis}[
            ybar,
            ymode = log,
            walltimestyle_proj,
            xlabel= {Size $n$} \\ {$M\dot{X}(t)M^T = A X(t) M^T +  M X(t) A^T + B B ^T,$} \\ {$X(0)=0,$}  \\  {$t\in [0,4500]$.},
            ]
            \addplot[fill=red,      discard if not={solver}{ode45}]     table[x={instance}, y={time}, col sep=comma] {results/timings_none.dat};
            \addplot[fill=orange,   discard if not={solver}{ode23}]     table[x={instance}, y={time}, col sep=comma] {results/timings_none.dat};
            \addplot[fill=magenta,  discard if not={solver}{ode113}]    table[x={instance}, y={time}, col sep=comma] {results/timings_none.dat};
            \addplot[fill=black,    discard if not={solver}{ode15s}]    table[x={instance}, y={time}, col sep=comma] {results/timings_none.dat};
            \addplot[fill=brown,    discard if not={solver}{ode23s}]    table[x={instance}, y={time}, col sep=comma] {results/timings_none.dat};
            \addplot[fill=blue,     discard if not={solver}{ode23t}]    table[x={instance}, y={time}, col sep=comma] {results/timings_none.dat};
            \addplot[fill=cyan,     discard if not={solver}{ode23tb}]   table[x={instance}, y={time}, col sep=comma] {results/timings_none.dat};
        \end{axis}
    \end{tikzpicture}
\end{minipage}
\hfill
\begin{minipage}{0.5\textwidth}
    \centering
    \begin{tikzpicture}[baseline]
        \begin{axis}[
            ybar,
            ymode = log,
            walltimestyle_proj,
            xlabel= {Size $n$} \\ {$M^T\dot{X}(t)M = A^T X(t) M +  M^T X(t) A + C^T C,$} \\ {$X(0)=0,$} \\ {$t\in [0,4500]$.},
            ]
            \addplot[fill=red,      discard if not={solver}{ode45}]     table[x={instance}, y={time}, col sep=comma] {results/timings_trans.dat};
            \addplot[fill=orange,   discard if not={solver}{ode23}]     table[x={instance}, y={time}, col sep=comma] {results/timings_trans.dat};
            \addplot[fill=magenta,  discard if not={solver}{ode113}]    table[x={instance}, y={time}, col sep=comma] {results/timings_trans.dat};
            \addplot[fill=black,    discard if not={solver}{ode15s}]    table[x={instance}, y={time}, col sep=comma] {results/timings_trans.dat};
            \addplot[fill=brown,    discard if not={solver}{ode23s}]    table[x={instance}, y={time}, col sep=comma] {results/timings_trans.dat};
            \addplot[fill=blue,     discard if not={solver}{ode23t}]    table[x={instance}, y={time}, col sep=comma] {results/timings_trans.dat};
            \addplot[fill=cyan,     discard if not={solver}{ode23tb}]   table[x={instance}, y={time}, col sep=comma] {results/timings_trans.dat};
        \end{axis}
    \end{tikzpicture}
\end{minipage}
\\[0.5em]
\begin{minipage}[c]{\textwidth}
    \centering
    \begin{tikzpicture}
        \begin{axis}[
                ybar,%
                hide axis,
                xmin=10,
                xmax=30,
                ymin=0,
                ymax=1.0,
                legend style={draw=white!20!black},
                legend columns =7,
                legend image post style={},
            ]
            \addlegendimage{fill=red};
            \addlegendentry{\texttt{ode45}};

            \addlegendimage{fill=orange};
            \addlegendentry{\texttt{ode23}};

            \addlegendimage{fill=magenta};
            \addlegendentry{\texttt{ode113}}

            \addlegendimage{fill=black};
            \addlegendentry{\texttt{ode15s}};

            \addlegendimage{fill=brown};
            \addlegendentry{\texttt{ode23s}}

            \addlegendimage{fill=blue};
            \addlegendentry{\texttt{ode23t}}

            \addlegendimage{fill=cyan};
            \addlegendentry{\texttt{ode23tb}}
        \end{axis}
    \end{tikzpicture}
    
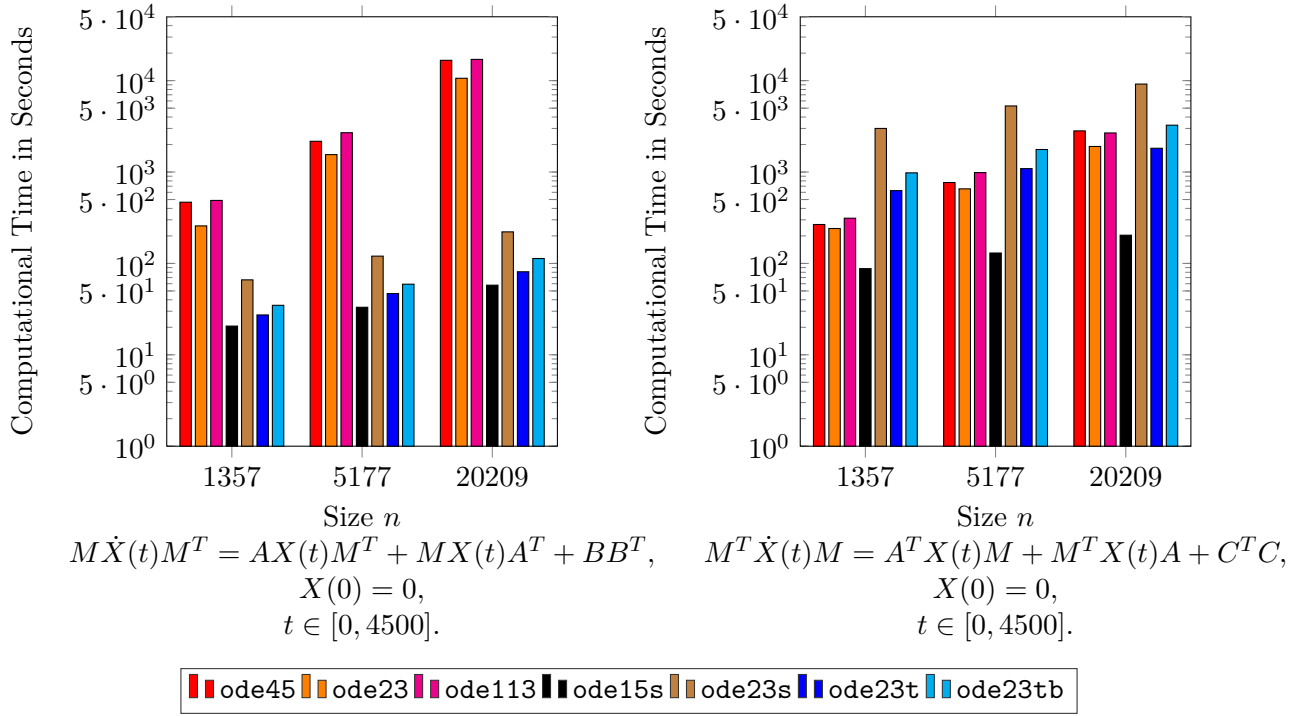
\captionof{figure}[]{Timings of the different \matlab{} ODE solvers for the solution of the projected time-dependent factors as defined in Equation~\eqref{eq:projected_td_factors}.}\label{fig:odesolvtimings}
\end{minipage}
\\[0.5em]
\begin{minipage}{0.5\textwidth}
    \centering
    \begin{tikzpicture}[baseline]
        \begin{axis}[
            ybar,
            ymode = log,
            walltimestyle_bdf,
            xlabel=Step Size $h$ \\ {$M\dot{X}(t)M^T = A X(t) M^T +  M X(t) A^T + B B ^T,$} \\ {$X(0)=0,$}  \\ {$t\in [0,100]$ and $n=1357$.},
            ]
            \addplot[fill=\BDFonecolor,     discard if not={solver}{BDF1}] table[ x={stepsize}, y={wtime}, col sep=comma] {results_ldlt/timings_none.dat};
            \addplot[fill=\BDFtwocolor,     discard if not={solver}{BDF2}] table[ x={stepsize}, y={wtime}, col sep=comma] {results_ldlt/timings_none.dat};
            \addplot[fill=\BDFthreecolor,   discard if not={solver}{BDF3}] table[ x={stepsize}, y={wtime}, col sep=comma] {results_ldlt/timings_none.dat};
            \addplot[fill=\BDFfourcolor,    discard if not={solver}{BDF4}] table[ x={stepsize}, y={wtime}, col sep=comma] {results_ldlt/timings_none.dat};
            \addplot[fill=\BDFfivecolor,    discard if not={solver}{BDF5}] table[ x={stepsize}, y={wtime}, col sep=comma] {results_ldlt/timings_none.dat};
            \addplot[fill=\BDFsixcolor,     discard if not={solver}{BDF6}] table[ x={stepsize}, y={wtime}, col sep=comma] {results_ldlt/timings_none.dat};
        \end{axis}
    \end{tikzpicture}
\end{minipage}
\hfill
\begin{minipage}{0.5\textwidth}
    \centering
    \begin{tikzpicture}[baseline]
        \begin{axis}[
            ybar,
            ymode = log,
            walltimestyle_bdf,
            xlabel=Step Size $h$ \\ {$M^T\dot{X}(t)M = A^T X(t) M +  M^T X(t) A + C^T C,$} \\ {$X(0)=0,$}  \\ {$t\in [0,100]$ and $n=1357$.},
            ]
            \addplot[fill=\BDFonecolor,     discard if not={solver}{BDF1}] table[ x={stepsize}, y={wtime}, col sep=comma] {results_ldlt/timings_trans.dat};
            \addplot[fill=\BDFtwocolor,     discard if not={solver}{BDF2}] table[ x={stepsize}, y={wtime}, col sep=comma] {results_ldlt/timings_trans.dat};
            \addplot[fill=\BDFthreecolor,   discard if not={solver}{BDF3}] table[ x={stepsize}, y={wtime}, col sep=comma] {results_ldlt/timings_trans.dat};
            \addplot[fill=\BDFfourcolor,    discard if not={solver}{BDF4}] table[ x={stepsize}, y={wtime}, col sep=comma] {results_ldlt/timings_trans.dat};
            \addplot[fill=\BDFfivecolor,    discard if not={solver}{BDF5}] table[ x={stepsize}, y={wtime}, col sep=comma] {results_ldlt/timings_trans.dat};
            \addplot[fill=\BDFsixcolor,     discard if not={solver}{BDF6}] table[ x={stepsize}, y={wtime}, col sep=comma] {results_ldlt/timings_trans.dat};
        \end{axis}
    \end{tikzpicture}
\end{minipage}
\\[0.5em]
\begin{minipage}{\textwidth}
    \begin{center}
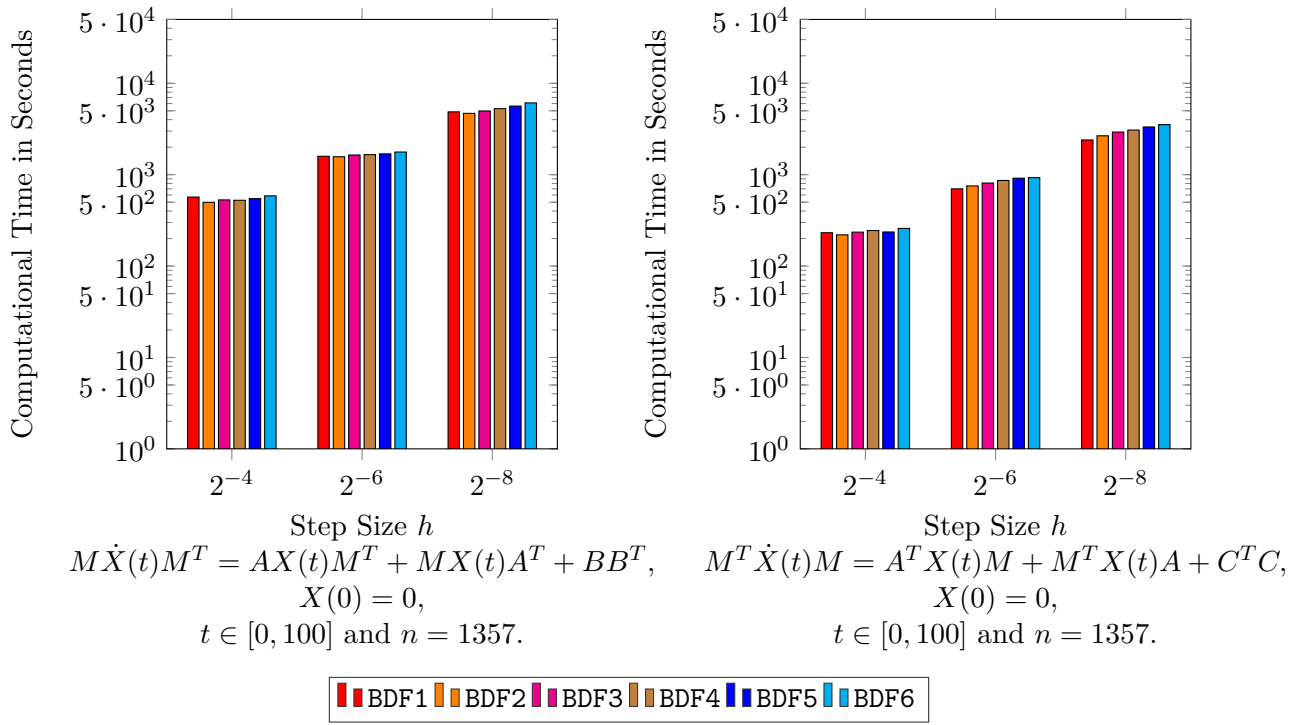

        \begin{tikzpicture}
            \begin{axis}[
                    ybar,%
                    hide axis,
                    xmin=10,
                    xmax=30,
                    ymin=0,
                    ymax=1.0,
                    legend style={draw=white!20!black},
                    legend columns =7,
                    legend image post style={},
                ]
                \addlegendimage{fill=\BDFonecolor};
                \addlegendentry{\texttt{BDF1}};

                \addlegendimage{fill=\BDFtwocolor};
                \addlegendentry{\texttt{BDF2}};

                \addlegendimage{fill=\BDFthreecolor};
                \addlegendentry{\texttt{BDF3}}

                \addlegendimage{fill=\BDFfourcolor};
                \addlegendentry{\texttt{BDF4}};

                \addlegendimage{fill=\BDFfivecolor};
                \addlegendentry{\texttt{BDF5}}

                \addlegendimage{fill=\BDFsixcolor};
                \addlegendentry{\texttt{BDF6}}

            \end{axis}
        \end{tikzpicture}
    \end{center}
    \captionof{figure}[]{Timings of the different BDF/ADI solvers.}\label{fig:ldltsolvtimings}
\end{minipage}

\subsection{Comparison with Backward Differentiation Formulas}\label{sec_numerical_res_comp_bdf}
To benchmark our method, we have run comparison tests with the \matlab{} implementation of the Backward Differentiation Formulas / Alternating Direction Implicit (BDF/ADI) scheme as developed in~\cite{LLanSS17} (see also the Appendix~\ref{sec_appendix_bdf} for a short summary of the algorithm).
The numerical experiments were conducted on the same machine, with the same \matlab{} version and the same model as described in Section~\ref{sec_numerical_res_setup}.

In contrast to the Algorithm~\ref{alg:projection} that needs to solve only a single algebraic Lyapunov equation for the initialization, the BDF/ADI approach solves a Lyapunov equation in every time step. Moreover, the numerical solution is stored in $LDL^T$-format, i.e., in terms of the factors of
$X(t_k) \approx L_k D_k L_k^T$ and $L_k \in \Rmat{n}{l_k}$, which grow at least linearly with the size of $n$ and the number of time steps.  For this reason, we had to restrict our numerical experiments with BDF/ADI to the interval $[0,100]$ and consider only the model of the smallest size $n=1357$.
As for the test of our Algorithm~\ref{alg:projection} in Appendix~\ref{sec_appendix_numerical_res_proj}, for computing the error, we used an exact solution $X_{ref}$ based on the spectral decomposition of $(A,M)$.

We have compared various BDF methods that we
abbreviated as \texttt{BDF1}, \texttt{BDF2}, \texttt{BDF3}, \texttt{BDF4},  \texttt{BDF5} and \texttt{BDF6}, where the number denotes the order $s$ of the method.
We used the constant time step sizes $\tau_k := h \in \{2^{-4}, 2^{-6}, 2^{-8}\}$ for our computations.
In Appendices~\ref{sec_appendix_numerical_res_bdf_dlec} and~\ref{sec_appendix_numerical_res_bdf_dleo}, we plot the relative and absolute errors compared to the solution
obtained by the spectral decomposition, c.f. Figures~\ref{fig:appendix:bdf:rel:none:k4}--\ref{fig:appendix:bdf:abs:none:k8} and~\ref{fig:appendix:bdf:rel:trans:k4}--\ref{fig:appendix:bdf:abs:trans:k8}.
. For comparison, we also plot the error of the numerical solution obtained by Algorithm~\ref{alg:projection} and the \matlab{} ODE solver
\texttt{ode15s}.
We list the computational times for performing the BDF/ADI method based computations in Section~\ref{sec_numerical_res_comptime}.
\\
As the actual solution $X$ converges towards the solution of the algebraic Lyapunov equation, also the numerical errors show a decay towards a certain error level. Decreasing the step size lowers this level accordingly. This effect is visible only for methods of lower order ($s\leq2$). For higher orders, the error level stagnates and the error rather shows oscillations, which are likely due to errors in the solution of the Lyapunov equations. Since all BDF methods are stiff solvers, there are no oscillations of higher frequency to observe. The error levels are comparable to the level reached by our approach with \texttt{ode15s}.

The computational timings for the BDF/ADI solvers are nearly the same for same step sizes; cf. Figure~\ref{fig:ldltsolvtimings}. Accordingly, the higher order methods clearly outperform the low-order methods.

As for the comparison to our approach, we note that the reported timings were for the longer time-interval $[0,4500]$. However, since the transient behavior of the solution is confined to a short initial phase and since the ODE solvers have a step size control, a restriction to a shorter interval would not change the timings by much.

For the test case with~\eqref{eq:numtestDLEc}, the BDF/ADI schemes achieved the same accuracy as our approach already with the coarsest step size; see Figures~\ref{fig:appendix:bdf:rel:none:k4}--\ref{fig:appendix:bdf:abs:none:k8}. Thus, we compare the timings in the left most columns in Figure~\ref{fig:odesolvtimings} ($n=1357$) and Figure~\ref{fig:ldltsolvtimings} ($h=2^{-4}$) to conclude that our approach with the best performing ODE solver is about $25$ times faster.

As for the test case with~\eqref{eq:numtestDLEo}, the BDF/ADI schemes reached the same accuracy only for the finest chosen step size; see Figures~\ref{fig:appendix:bdf:rel:trans:k4}--\ref{fig:appendix:bdf:abs:trans:k8}. Thus, comparing the right most column in Figure~\ref{fig:ldltsolvtimings} ($h=2^{-8}$) 
to the left column in the right plot of Figure~\ref{fig:odesolvtimings} ($n=1357$), 
we find that our Algorithm~\ref{alg:projection}, again, is faster by a factor of $30$.


To conclude the comparison, we add the following remarks. With the costs of solving Lyapunov equations, apart from the increased memory requirements,
the BDF/ADI scheme will also become significantly more costly for larger system sizes. As opposed to our Algorithm~\ref{alg:projection} however, the BDF/ADI scheme also applies for the time varying case as well as for nonzero initial conditions. Moreover, as \texttt{ode15s} in fact uses the BDF formulas but with variable order $s$ and variable time step, one may consider integrating similar error control mechanisms in the BDF/ADI approach to improve performance.

\section{Conclusions}\label{sec_conclusions}
We presented several solution formulas for the differential Sylvester and Lyapunov equations. For the autonomous stable differential Lyapunov equation, we proposed a numerical algorithm that combined certain aspects derived from the solution representations. The main feature of the algorithm is the projection of the time dependent part onto a suitable subspace of lower dimension. Only this makes the numerical solution feasible in terms of memory requirements. As for the computational time, the projected system can be directly solved with optimized ODE solvers like the \texttt{ode-suite} in \matlab{}. Moreover, its structure allows for column-wise computation of the factors and, thus, for straight-forward parallelization.

We illustrated the performance of the algorithm in an example that was derived from a finite-element discretization of a heat equation. The achieved accuracy is fully satisfactory. The greatest benefit, also in contrast to existing numerical schemes, is the low memory requirement.

The possible extension to the unstable or unsymmetric as well as to the non-autonomous case is not straightforward since the space in which the solution evolves has not been found to span a suitable invariant Krylov subspace and is possibly of high dimension.
Moreover, the solution of $AX=C$ (which is the special case of $AX+XB=C$ with $B=0$) is unlikely to span an $A$-invariant subspace at all, meaning that the span of the solution of the algebraic Sylvester equation is not suited for the differential equation. Thus, for the general case, a different ansatz is needed. The same is true for time-dependent coefficients or non-zero initial conditions. Here, a remedy might be structured approaches for the case that time dependence comes, e.g., from a low-rank update.

In the unstable case, apart from the fact that the algebraic Lyapunov equation may not have a unique solution, there can be modes that grow exponentially in time. For this case it may be worth investigating whether a projector that identifies the stable part of the underlying mathematical model can be efficiently incorporated in the proposed solution approach.


\pagebreak


\appendix
\section{Numerical Results Projection Approach}\label{sec_appendix_numerical_res_proj}

\newcommand{\ErrorPlotsProjScale}{0.9\textwidth}
\pgfplotsset{
    axisrellongtimestyle/.style         = {align = center, xscale = 1.0, yscale = 1.0, width = \ErrorPlotsProjScale, yminorticks = false, ylabel = {rel. 2-norm error},    xlabel = {$t\in [0,4500]$}},
    axisabslongtimestyle/.style         = {align = center, xscale = 1.0, yscale = 1.0, width = \ErrorPlotsProjScale, yminorticks = false, ylabel = {abs. 2-norm error},    xlabel = {$t\in [0,4500]$}},
    axisrelshorttimestyle/.style        = {align = center, xscale = 1.0, yscale = 1.0, width = \ErrorPlotsProjScale, yminorticks = false, ylabel = {rel. 2-norm error},    xlabel = {$t\in [0,10]$}},
    axisabsshorttimestyle/.style        = {align = center, xscale = 1.0, yscale = 1.0, width = \ErrorPlotsProjScale, yminorticks = false, ylabel = {abs. 2-norm error},    xlabel = {$t\in [0,10]$}},
    axis2nrmsolshorttimestyle/.style    = {align = center, xscale = 1.0, yscale = 1.0, width = \ErrorPlotsProjScale, yminorticks = false, ylabel = {2-norm solution},      xlabel = {$t\in [0,10]$}},
    axis2nrmsollongtimestyle/.style     = {align = center, xscale = 1.0, yscale = 1.0, width = \ErrorPlotsProjScale, yminorticks = false, ylabel = {2-norm solution},      xlabel = {$t\in [0,4500]$}},
    proj_line_style/.style              = {thick},
    proj_rel_yaxis_none/.style          = {ymin = 1e-9,     ymax = 1e-4,    ytickten={-9,-8,...,-4}},
    proj_abs_yaxis_none/.style          = {ymin = 1e-13,    ymax = 1e-9,    ytickten={-13,-12,...,-9}},
    proj_rel_yaxis_trans/.style         = {ymin = 1e-13,    ymax = 1e-8,    ytickten={-13,-12,...,-8}},
    proj_abs_yaxis_trans/.style         = {ymin = 1e-3,     ymax = 1e+2,    ytickten={-3,-2,...,2}},
    true_sol_abs_trans/.style           = {ymin = 1e+8,     ymax = 5e+11,   ytickten={8.0,8.5,...,11.5}},
    true_sol_abs_none/.style            = {ymin = 1e-6,     ymax = 1e-2,    ytickten={-6.0,-5.5,...,-2.0}},
}

\newcommand{\ErrorPlotsProjLegend}[1]{
    \begin{minipage}[#1]{\textwidth}
        \centering
        \begin{tikzpicture}
            \begin{axis}[%
                    hide axis,
                    xmin=10,
                    xmax=30,
                    ymin=0,
                    ymax=1.0,
                    line width=0.1mm,
                    legend style={draw=white!20!black},
                    legend columns =7,
                    legend image post style={line width=0.5mm},
                ]
                \addlegendimage{red};
                \addlegendentry{\texttt{ode45}};

                \addlegendimage{orange};
                \addlegendentry{\texttt{ode23}};

                \addlegendimage{magenta};
                \addlegendentry{\texttt{ode113}}

                \addlegendimage{black};
                \addlegendentry{\texttt{ode15s}};

                \addlegendimage{brown};
                \addlegendentry{\texttt{ode23s}}

                \addlegendimage{blue};
                \addlegendentry{\texttt{ode23t}}

                \addlegendimage{cyan};
                \addlegendentry{\texttt{ode23tb}}
            \end{axis}
        \end{tikzpicture}
    \end{minipage}
}

\newcommand{\ErrorPlotsProjMiniPageWidth}{0.45\textwidth}

\newcommand{\ErrorPlotsProjLong}[2]{
    \begin{minipage}{\ErrorPlotsProjMiniPageWidth}
        \captionsetup{type=figure}
        \centering
        \begin{tikzpicture}[baseline]
            \begin{semilogyaxis}[axisrellongtimestyle, proj_rel_yaxis_#2]
                \pgfplotsforeachungrouped \forcol/\forode in {red/ode45,orange/ode23,magenta/ode113,black/ode15s,brown/ode23s,blue/ode23t,cyan/ode23tb}{
                        \expandafter\addplot\expandafter[\forcol, proj_line_style]
                        table [skip first n=15, x={Time}, y={RELNRM2ERR},  col sep=comma, header=true] {results/rail#1_\forode_op_#2_false_local_parallel_FOURIER_false_SPECTRALDECOMP_true.dat};
                    }
            \end{semilogyaxis}
        \end{tikzpicture}
        \captionof{figure}{Relative 2-norm error of the approximation.}\label{fig:appendix:proj:rel:long:#1:#2}
    \end{minipage}
    \hfill
    \begin{minipage}{\ErrorPlotsProjMiniPageWidth}
        \captionsetup{type=figure}
        \centering
        \begin{tikzpicture}[baseline]
            \begin{semilogyaxis}[axisabslongtimestyle, proj_abs_yaxis_#2]
                \pgfplotsforeachungrouped \forcol/\forode in {red/ode45,orange/ode23,magenta/ode113,black/ode15s,brown/ode23s,blue/ode23t,cyan/ode23tb}{
                        \expandafter\addplot\expandafter[\forcol, proj_line_style]
                        table [skip first n=15, x={Time}, y={ABSNRM2ERR},  col sep=comma, header=true] {results/rail#1_\forode_op_#2_false_local_parallel_FOURIER_false_SPECTRALDECOMP_true.dat};
                    }
            \end{semilogyaxis}
        \end{tikzpicture}
        \captionof{figure}{Absolute 2-norm error of the approximation.}\label{fig:appendix:proj:abs:long:#1:#2}
    \end{minipage}
}

\newcommand{\ErrorPlotsProjShort}[2]{
    \begin{minipage}{\ErrorPlotsProjMiniPageWidth}
        \captionsetup{type=figure}
        \centering
        \begin{tikzpicture}[baseline]
            \begin{semilogyaxis}[axisrelshorttimestyle, proj_rel_yaxis_#2]
                \pgfplotsforeachungrouped \forcol/\forode in {red/ode45,orange/ode23,magenta/ode113,black/ode15s,brown/ode23s,blue/ode23t,cyan/ode23tb}{
                        \expandafter\addplot\expandafter[\forcol, proj_line_style]
                        table [skip first n=15, x={Time}, y={RELNRM2ERR},  col sep=comma, header=true, restrict x to domain=0:10] {results/rail#1_\forode_op_#2_false_local_parallel_FOURIER_false_SPECTRALDECOMP_true.dat};
                    }
            \end{semilogyaxis}
        \end{tikzpicture}
        \captionof{figure}{Relative 2-norm error of the approximation.}\label{fig:appendix:proj:rel:short:#1:#2}
    \end{minipage}
    \hfill
    \begin{minipage}{\ErrorPlotsProjMiniPageWidth}
        \captionsetup{type=figure}
        \centering
        \begin{tikzpicture}[baseline]
            \begin{semilogyaxis}[axisabsshorttimestyle, proj_abs_yaxis_#2]
                \pgfplotsforeachungrouped \forcol/\forode in {red/ode45,orange/ode23,magenta/ode113,black/ode15s,brown/ode23s,blue/ode23t,cyan/ode23tb}{
                        \expandafter\addplot\expandafter[\forcol, proj_line_style]
                        table [skip first n=15, x={Time}, y={ABSNRM2ERR},  col sep=comma, header=true, restrict x to domain=0:10] {results/rail#1_\forode_op_#2_false_local_parallel_FOURIER_false_SPECTRALDECOMP_true.dat};
                    }
            \end{semilogyaxis}
        \end{tikzpicture}
        \captionof{figure}{Absolute 2-norm error of the approximation.}\label{fig:appendix:proj:abs:short:#1:#2}
    \end{minipage}
}

\newcommand{\ErrorPlotsSolution}[2]{
    \begin{minipage}{\ErrorPlotsProjMiniPageWidth}
        \captionsetup{type=figure}
        \centering
        \begin{tikzpicture}[baseline]
            \begin{semilogyaxis}[axis2nrmsollongtimestyle, true_sol_abs_#2]
                \addplot [color=green!75!black, proj_line_style]
                table [skip first n=15, x={Time}, y={NRM2XTRUE}, col sep=comma, header=true] {results/rail#1_ode45_op_#2_false_local_parallel_FOURIER_false_SPECTRALDECOMP_true.dat};
            \end{semilogyaxis}
        \end{tikzpicture}
        \captionof{figure}{2-norm of the reference solution.}\label{fig:appendix:proj:norm:long:#1:#2}
    \end{minipage}
    \hfill
    \begin{minipage}{\ErrorPlotsProjMiniPageWidth}
        \captionsetup{type=figure}
        \centering
        \begin{tikzpicture}[baseline]
            \begin{semilogyaxis}[axis2nrmsolshorttimestyle, true_sol_abs_#2]
                \addplot [color=green!75!black, proj_line_style]
                table [skip first n=15, x={Time}, y={NRM2XTRUE}, col sep=comma, header=true, restrict x to domain=0:10] {results/rail#1_ode45_op_#2_false_local_parallel_FOURIER_false_SPECTRALDECOMP_true.dat};
            \end{semilogyaxis}
        \end{tikzpicture}
        \captionof{figure}{2-norm of the reference solution.}\label{fig:appendix:proj:norm:short:#1:#2}
    \end{minipage}
}

\newcommand{\ErrorPlotsProjGroup}[2]{
    \ErrorPlotsProjLong{#1}{#2}
    \\[0.5em]
    \ErrorPlotsProjShort{#1}{#2}
    \\[0.5em]
    \ErrorPlotsProjLegend{c}
    \\[0.5em]
    \ErrorPlotsSolution{#1}{#2}
}

\subsection{Results for the Differential Lyapunov Equation}\label{sec_appendix_numerical_res_proj_dlec}
\begin{center}$n=1357$ and $M\dot{X}(t)M^T = A X(t) M^T +  M X(t) A^T + B B ^T,\ X(0)=0$. \end{center}
\ErrorPlotsProjGroup{1357}{none}
\pagebreak

\begin{center}$n=5177$ and $M\dot{X}(t)M^T = A X(t) M^T +  M X(t) A^T + B B ^T,\ X(0)=0$. \end{center}
\ErrorPlotsProjGroup{5177}{none}
\pagebreak

\begin{center}$n=20209$ and $M\dot{X}(t)M^T = A X(t) M^T +  M X(t) A^T + B B ^T,\ X(0)=0$. \end{center}
\ErrorPlotsProjGroup{20209}{none}
\pagebreak

\subsection{Results for the Transposed Differential Lyapunov Equation}\label{sec_appendix_numerical_res_proj_dleo}
\begin{center}$n=1357$ and $M^T\dot{X}(t)M = A^T X(t) M +  M^T X(t) A + C^T C,\ X(0)=0$. \end{center}
\ErrorPlotsProjGroup{1357}{trans}
\pagebreak

\begin{center}$n=5177$ and $M^T\dot{X}(t)M = A^T X(t) M +  M^T X(t) A + C^T C,\ X(0)=0$. \end{center}
\ErrorPlotsProjGroup{5177}{trans}
\pagebreak

\begin{center}$n=20209$ and $M^T\dot{X}(t)M = A^T X(t) M +  M^T X(t) A + C^T C,\ X(0)=0$. \end{center}
\ErrorPlotsProjGroup{20209}{trans}
\pagebreak

\section{Backward Differentiation Formulas}\label{sec_appendix_bdf}
\newcommand{\ErrorPlotsBDFScale}{0.9\textwidth}
\pgfplotsset{
    bdf_rel_style/.style        = {align = center, xscale = 1.0, yscale = 1.0, width = \ErrorPlotsBDFScale, yminorticks = false, ylabel = {rel. 2-norm error}},
    bdf_abs_style/.style        = {align = center, xscale = 1.0, yscale = 1.0, width = \ErrorPlotsBDFScale, yminorticks = false, ylabel = {abs. 2-norm error}},
    bdf_line_style/.style       = {thick},
    bdf_rel_yaxis_none/.style   = {ymin=1e-10,  ymax = 1e-1,    ytickten = {-10,-8,...,-1}},
    bdf_abs_yaxis_none/.style   = {ymin=1e-18,  ymax = 1e-5,    ytickten = {-18,-16,...,-5}},
    bdf_rel_yaxis_trans/.style  = {ymin=1e-12,  ymax = 1e-1,    ytickten = {-12,-10,...,-1}},
    bdf_abs_yaxis_trans/.style  = {ymin=1e-6,   ymax = 1e8,     ytickten = {-6,-4,...,8}},
}

\newcommand{\ErrorPlotsBDFLegend}[1]{
    \begin{minipage}[#1]{\textwidth}
        \centering
        \begin{tikzpicture}
            \begin{axis}[%
                    hide axis,
                    xmin=10,
                    xmax=30,
                    ymin=0,
                    ymax=1.0,
                    line width=0.1mm,
                    legend style={draw=white!20!black},
                    legend columns =8,
                    legend image post style={line width=0.5mm},
                ]
                \addlegendimage{black};
                \addlegendentry{\texttt{ode15s}};

                \addlegendimage{\BDFonecolor};
                \addlegendentry{\texttt{BDF1}};

                \addlegendimage{\BDFtwocolor};
                \addlegendentry{\texttt{BDF2}}

                \addlegendimage{\BDFthreecolor};
                \addlegendentry{\texttt{BDF3}};

                \addlegendimage{\BDFfourcolor};
                \addlegendentry{\texttt{BDF4}}

                \addlegendimage{\BDFfivecolor};
                \addlegendentry{\texttt{BDF5}}

                \addlegendimage{\BDFsixcolor};
                \addlegendentry{\texttt{BDF6}}
            \end{axis}
        \end{tikzpicture}
    \end{minipage}
}

\newcommand{\ErrorPlotsBDFjMiniPageWidth}{0.45\textwidth}

\newcommand{\ErrorPlotsBDF}[2]{
    \begin{minipage}{\ErrorPlotsBDFjMiniPageWidth}
        \captionsetup{type=figure}
        \centering
        \begin{tikzpicture}[baseline]
            \begin{semilogyaxis}[bdf_rel_style, bdf_rel_yaxis_#1, xlabel = {$t\in [0,100],\ h=2^{-#2}$}]
                \pgfplotsforeachungrouped \forcol/\forode in {\BDFonecolor/BDF1,\BDFtwocolor/BDF2,\BDFthreecolor/BDF3,\BDFfourcolor/BDF4,\BDFfivecolor/BDF5,\BDFsixcolor/BDF6}{
                        \expandafter\addplot\expandafter[\forcol, bdf_line_style]
                        table [mark=none, skip first n=9, each nth point=1, x={Time}, y={RELNRM2ERR},  col sep=comma, header=true]  {results_ldlt/rail1357_\forode_T100.000000_k#2_op_#1_sparse.dat};
                    }
                \addplot[black, bdf_line_style]
                table [skip coords between index={0}{0}, skip first n=15, x={Time}, y={RELNRM2ERR},  col sep=comma, header=true] {results_ldlt/rail1357_ode15s_op_#1_ldlt_comp.dat};
            \end{semilogyaxis}
        \end{tikzpicture}
        \captionof{figure}{Relative 2-norm error of the {BDF/ADI} approximation.}\label{fig:appendix:bdf:rel:#1:k#2}
    \end{minipage}
    \hfill
    \begin{minipage}{\ErrorPlotsBDFjMiniPageWidth}
        \captionsetup{type=figure}
        \centering
        \begin{tikzpicture}[baseline]
            \begin{semilogyaxis}[bdf_abs_style, bdf_abs_yaxis_#1, xlabel = {$t\in [0,100],\ h=2^{-#2}$}]
                \pgfplotsforeachungrouped \forcol/\forode in {\BDFonecolor/BDF1,\BDFtwocolor/BDF2,\BDFthreecolor/BDF3,\BDFfourcolor/BDF4,\BDFfivecolor/BDF5,\BDFsixcolor/BDF6}{
                        \expandafter\addplot\expandafter[\forcol, bdf_line_style]
                        table [mark=none, skip first n=9, each nth point=1, x={Time}, y={ABSNRM2ERR},  col sep=comma, header=true]  {results_ldlt/rail1357_\forode_T100.000000_k#2_op_#1_sparse.dat};
                    }
                \addplot[black, bdf_line_style]
                table [skip coords between index={0}{0}, skip first n=15, x={Time}, y={ABSNRM2ERR},  col sep=comma, header=true] {results_ldlt/rail1357_ode15s_op_#1_ldlt_comp.dat};
            \end{semilogyaxis}
        \end{tikzpicture}
        \captionof{figure}{Absolute 2-norm error of the {BDF/ADI} approximation.}\label{fig:appendix:bdf:abs:#1:k#2}
    \end{minipage}
}

\newcommand{\ErrorPlotsBDFGroup}[4]{
    \ErrorPlotsBDF{#1}{#2}
    \\[0.5em]
    \ErrorPlotsBDF{#1}{#3}
    \\[0.5em]
    \ErrorPlotsBDF{#1}{#4}
    \\[0.5em]
    \ErrorPlotsBDFLegend{c}
}

We consider Backward Differentiation Formulas (BDF) for differential Lyapunov equations~\cite{morLanSS16,morLanSS15}.
Let $0=t_0<t_1<\cdots<t_{N} =T$ be a decomposition of the interval $[0,T]$. We define the step size $\tau_k = t_{k} - t_{k-1}$ for $k=1,\ldots,N$. \\
The $s$-step BDF method applied to the DLE~\ref{eqn:autodle} is given by
\begin{align*}
    \sum\limits_{j=0}^{s} \alpha_j X_{k-j} = \tau_k \beta \left(A^T X_k + X_k A + B{B^T} \right),
\end{align*}
where $\alpha_j$ and $\beta$ are coefficients of the BDF method and can be found in~\cite{HaiNW87}.
The parameter $s$ is the order of the BDF method. We recall that for $s>6$, the method is not numerical stable, and for $s=1$, the BDF method coincides with the implicit Euler method.
A minor rearrangement shows that the current iterate $X_k$ can be obtained as the solution of the algebraic Lyapunov equation
\begin{align}\label{eq:bdf-lyap-eqn}
    {\left( \tau_k \beta A - \frac{\alpha_0}{2} E_{n,n} \right)}^T X_k + X_k \left( \tau_k \beta A - \frac{\alpha_0}{2} E_{n,n}\right) = -\tau_k \beta BB^T + \sum\limits_{j=1}^{s} \alpha_j X_{k-j}.
\end{align}
Since for $s\geq 2$, certain coefficients $\alpha_j$, $j\geq 1$ are positive, the algebraic Lyapunov equation~\eqref{eq:bdf-lyap-eqn} has a symmetric but possibly indefinite right-hand side, which makes the standard ADI method infeasible.
For this reason a $LDL^T$-decomposition for the numerical solution is proposed and suitable modifications of the ADI method have been developed;~\cite{LLanMS15,BenLT09}: Assume that $X_i \approx L_i D_i L_i^T$ for $i=0,\ldots, k-1$, $L_i\in \Rmat{n}{l_i},~ D_i \in \Rmat{l_i}{l_i}$ and $l_i \ll n$, then the right hand side can be factored as
\begin{align*}
    -\tau_k \beta BB^T + \sum\limits_{j=1}^{s} \alpha_j X_{k-j} & \approx  -G_k S_k G_k^T,     \\
    G_k                                                         & = \begin{bmatrix} B, L_{k-1}, \ldots, L_{k-s} \end{bmatrix}, \\
    S_k                                                         & =
    \begin{bmatrix}
        \tau_k \beta E_{p,p} &                   &        &                   \\
                             & -\alpha_1 D_{k-1} &        &                   \\
                             &                   & \ddots &                   \\
                             &                   &        & -\alpha_s D_{k-s}
    \end{bmatrix}.
\end{align*}
Now the $LDL^T$-type ADI method can be used to determine $X_k \approx  L_k D_k L_k^T$.
The BDF/ADI methods can also be extended to generalized differential Lyapunov equations in a similar way~\cite{morLanSS16}.

\pagebreak
\subsection{Results for the Differential Lyapunov Equation}\label{sec_appendix_numerical_res_bdf_dlec}
\begin{center}$n=1357$ and $M\dot{X}(t)M^T = A X(t) M^T +  M X(t) A^T + B B ^T,\ X(0)=0$. \end{center}
\ErrorPlotsBDFGroup{none}{4}{6}{8}

\subsection{Results for the Transposed Differential Lyapunov Equation}\label{sec_appendix_numerical_res_bdf_dleo}
\begin{center}$n=1357$ and $M^T\dot{X}(t)M = A^T X(t) M +  M^T X(t) A + C^T C,\ X(0)=0$. \end{center}
\ErrorPlotsBDFGroup{trans}{4}{6}{8}


\bibliography{mybib/local,csc-bibfiles/csc,csc-bibfiles/mor,csc-bibfiles/software}

\end{document}